\documentclass[11pt,epic,eepic]{article}
\usepackage{fullpage, graphicx, latexsym, amssymb, amscd, psfrag}
\input amssym.def

\newcommand{\B}{{\mathbb B}}
\newcommand{\I}{{\mathbb I}}

\def\R{{{\rm I}\kern-.16em {\rm R}}}
\def\Z{{{\rm Z}\kern-.28em{\rm Z}}}
\def\N{{{\rm I}\kern-.16em {\rm N}}}
\def\C{{\rm C\kern-.48em\vrule width.06em height.57em depth-.02em \kern.48em}}

\newcommand{\qed}{\hfill $\Box$ \hfill \\}

\newcommand{\since}{.^{\kern-6.5pt{\textstyle\cdot\kern5pt\cdot}}~~}

\newcommand{\backs}{\backslash}

\def\span{{\rm span}}
\def\Exp{{\rm Exp}}
\def\Ideal{{\rm Ideal}}

\def\inpro#1{\langle#1\rangle}
\def\Set#1{\left\{\,#1\,\right\}}
\def\gp#1{\left(#1\right)}

\def\bd{\begin{description}}
\def\ed{\end{description}}

\def\be{\begin{equation}} 
\def\ee{\end{equation}}

\def\alp{\alpha}                

                \def\Gam{\Gamma}
\def\del{\delta}                
\def\eps{\varepsilon}

\def\lam{\lambda}               \def\Lam{\Lambda}
\def\sig{\sigma}

\def\calD{{\cal D}}

\def\calF{{\cal F}}

\def\calH{{\cal H}}
\def\calI{{\cal I}}
\def\calJ{{\cal J}}

\def\calP{{\cal P}}
\def\calQ{{\cal Q}}

\def\calZ{{\cal Z}}

\def\tilB{{\widetilde B}}

\def\tilQ{{\widetilde Q}}

\def\til#1{{\widetilde #1}}

\newenvironment{@abssec}[1]{%
     \if@twocolumn
       \section*{#1}%
     \else
       \vspace{.05in}\footnotesize
       \parindent .2in
         {\bfseries #1. }\ignorespaces
     \fi}
     {\if@twocolumn\else\par\vspace{.1in}\fi}
\newenvironment{keywords}{\begin{@abssec}{Key words}}{\end{@abssec}}
\newenvironment{AMS}{\begin{@abssec}{AMS subject classification}}{\end{@abssec}}

\newtheorem{theorem}{Theorem}[section] 
\newtheorem{proposition}[theorem]{Proposition}
\newtheorem{conjecture}[theorem]{Conjecture}
\newtheorem{lemma}[theorem]{Lemma}
\newtheorem{corollary}[theorem]{Corollary}
\newtheorem{result}[theorem]{Result}

\newtheorem{example}[theorem]{Example}
\newtheorem{definition}[theorem]{Definition}

\def\calQ{{\cal P}}

\def\implies{\Longrightarrow}
\def\le{\leq}

\def\calD{{\cal D}}

\def\calI{{\cal I}}
\def\calJ{{\cal J}}
\def\calP{{\cal P}}
\def\calZ{{\cal Z}}
\def\calF{{\cal F}}

\def\lam{\lambda}
\def\Lam{\Lambda}
\def\eps{\epsilon}

\def\alp{\alpha}

\def\codim{\mathop{\rm codim}\nolimits}

\def\eqbd{\mathop{{:}{=}}}
\def\bdeq{\mathop{{=}{:}}}
\def\least{_\downarrow}
\def\most{_\uparrow}
\def\bfH{{\bf H}}

\def\openC{{\rm C\kern-.48em\vrule width.06em height.6em depth-.02em 
                 \kern.48em}}
\def\openR{{{\rm I}\kern-.16em {\rm R}}}
\def\openZ{{{\rm Z}\kern-.28em{\rm Z}}}
\def\sZZ{{{\scriptstyle\rm Z}\kern-.24em{\scriptstyle\rm Z}}}
\def\openT{{{\rm T}\kern-.42em {\rm T}}}
\def\openH{{{\rm I}\kern-.16em {\rm H}}}
\def\openK{{{\rm I}\kern-.16em {\rm K}}}
\def\openL{{{\rm I}\kern-.16em {\rm L}}}
\def\openM{{{\rm I}\kern-.16em {\rm M}}}
\def\openN{{{\rm I}\kern-.16em {\rm N}}}
\def\openP{{{\rm I}\kern-.16em {\rm P}}}
\def\openB{{{\rm I}\kern-.16em {\rm B}}}
\def\oB{\openB}

\def\eqbd{\mathop{{:}{=}}}
\def\bdeq{\mathop{{=}{:}}}

\def\ee{{\rm e}}
\let\C\openC
\def\sC{{\rm C\kern-.38em\vrule width.06em height.45em depth-.02em 
                 \kern.3em}}
\let\N\openN

\let\R\openR
\let\Z\openZ
\let\B\openB
\def\eop{\hfill
        {\ \vbox{\hrule\hbox{\vrule height1.3ex\hskip0.8ex\vrule}\hrule}}
        \vskip 0.3cm \par}
\def\var{\mathop{\rm Var}\nolimits}

\def\spam{\mathop{\rm span}\nolimits}
\def\spa{\mathop{\rm span}\nolimits}

\def\rank{\mathop{\rm rank}\nolimits}
\def\ran{\mathop{\rm ran}\nolimits}

\def\Ideal{\mathop{\rm Ideal}\nolimits}
\def\ex{\mathop{\rm ex}\nolimits}
\def\val{\mathop{\rm val}\nolimits}
\def\belowrightarrow#1{{{{}\over\ #1\ }\kern-1.1em\to}}
\def\l2{{L_2}}

\def\bks{\backslash}

\def\Set#1{\{#1\}}
\def\nt{\noindent}
\def\Rn{\R^n}
\def\vol{{\rm vol}}
\def\Zn{\Z^n}
\def\bi{\begin{itemize}}
\def\ei{\end{itemize}}
\def\pzero{{\vphantom 0}}
\def\bprec{\mathop{\hbox{\bf{\LARGE$\prec$}}}\nolimits}
\def\bsucc{\mathop{\hbox{\bf{\LARGE$\succ$}}}\nolimits}
\def\subbsucc{\mathop{\hbox{\bf{\large$\succ$}}}\nolimits}

\begin{document}

\title{Zonotopal algebra}
\author{ Olga Holtz\thanks{Supported by the Sofja Kovalevskaja Research 
Prize of Alexander von Humboldt Foundation.} \\ Departments of Mathematics \\
University of California-Berkeley \\
\& Technische Universit\"at Berlin  \and Amos Ron\thanks{Supported by the US National Science 
Foundation under Grants ANI-0085984 and DMS-0602837,
by the National Institute of General Medical Sciences under Grant
NIH-1-R01-GM072000-01, and by the Vilas Foundation at the University of
Wisconsin.} \\
Department of Mathematics \& \\
Department of Computer Science \\
University of Wisconsin-Madison }
\date{January 24, 2008; updated February 16, 2011}
\maketitle

\begin{keywords} 
Multivariate polynomials, polynomial ideals, duality, grading, kernels of
differential operators, polynomial interpolation, box splines, zonotopes, 
hyperplane arrangements, matroids, graphs, parking functions,  
Tutte polynomial, Ehrhart polynomial,  Hilbert series. 
\end{keywords}

\begin{AMS} 
13F20, 13A02, 16W50, 16W60, 47F05, 47L20, 05B20, 05B35, 05B45, 05C50, 
52B05, 52B12, 52B20, 52C07, 52C35, 41A15, 41A63.
\end{AMS}

\begin{abstract} 
A wealth of geometric and combinatorial properties of a given linear 
endomorphism $X$ of $\R^N$ is captured in the study of its associated 
zonotope $Z(X)$, and, by duality, its associated hyperplane arrangement 
$\calH(X)$. This well-known line of study is particularly  interesting in case 
$n\eqbd\rank X \ll N$. We enhance this study to an algebraic level, and associate 
$X$ with three algebraic structures, referred herein as {\it external, 
central, and internal.\/} Each algebraic structure is given in terms of  a 
pair of homogeneous polynomial ideals in $n$ variables that are dual to 
each other: one encodes properties of the arrangement $\calH(X)$, while the 
other encodes by duality properties of the zonotope $Z(X)$. The algebraic 
structures are defined purely in terms of the combinatorial structure of
$X$, but are subsequently proved to be equally obtainable by applying suitable 
algebro-analytic operations to either of $Z(X)$ or $\calH(X)$.  The theory 
is universal in the sense that it requires no assumptions on the map $X$ 
(the only exception being that the algebro-analytic operations on $Z(X)$ 
yield sought-for results only in case $X$ is unimodular), and  provides
new tools that can be used in enumerative combinatorics, graph theory,
representation theory, polytope geometry, and approximation theory.  
\end{abstract}

\section{\label{sec_intro}Introduction}
\subsection{\label{sec_general}General}                                                                       

We are interested in combinatorial, geometric, algebraic and analytic
properties of low rank linear endomorphisms $X$ of $\R^N$.  This setup 
is relevant in quite a few areas in mathematics from linear algebra to
algebraic graph theory to semi-simple group representations to
approximation theory (box splines), and underlies interesting connections
among rather different mathematical problems.

Consider $X$ as a map from $\R^N$ to $\R^n$, and identify it with the 
columns of its matrix representation.  
Important geometric information about $X$ is captured by the image 
$$Z(X)\eqbd \{\sum_{x\in X}t_x x : \quad t\in [0,1]^X\}$$ 
of the unit cube $[0,1]^X$ under the action of $X$.  The resulting polytope 
is known as a {\it zonotope.\/} Zonotopes exhibit special symmetries that 
general polytopes lack. Underlying those special features is the fact that 
their normal cone fan is linear, i.e., is a (central) hyperplane 
arrangement. The duality between zonotopes and hyperplane arrangements is 
rich, and includes intriguing connections between the different tilings of 
the zonotope into sub-zonotopes on the one hand, and  the geometries obtained 
by translating the hyperplanes in the hyperplane arrangement on the other 
hand (see~\cite{McM_M}, \cite{McM_TAMS}, \cite{Shep_M}, \cite{Shep_CJM}, 
\cite[Chapter~7]{Ziegler_book}, \cite[Chapter~2]{BLVSWZ_book}, 
\cite{RGZ_CM}, \cite{S_notes}). While we briefly touch in 
Section~\ref{sec_prelim} on these known connections, the focus of this paper 
is neither on the linear algebra surrounding the map $X$, nor on the 
geometry and combinatorics of the zonotope $Z(X)$ {\it per se.}

The theory of {\it zonotopal algebra\/} that is developed in the current
article is algebraic. At its core one finds three pairs of
zero-dimensional homogeneous polynomial ideals in $n$ variables:
an external pair $(\calI_+(X),\calJ_+(X))$, a central pair
$(\calI(X),\calJ(X))$, and an internal pair $(\calI_-(X),\calJ_-(X))$.
The ideals within each pair are dual to each other; in particular, their
Hilbert series are identical. To keep this introduction brief, we do not
describe in depth the actual ingredients of the theory that is developed
here.  Instead, we present a number of results that capture the flavor of 
the general theory and its many potential applications.

The definition of the $\calI$-ideals
goes as follows. First, given $y\in \Rn$, let $p_y$ be the linear
form
$$p_y:\Rn\to\R\ : t\mapsto y\cdot t\eqbd \sum_{i=1}^ny(i)t(i).$$
Further, let
$$\calF(X)$$
be the set of {\it facet hyperplanes\/} of $X$, viz., $H\in \calF(X)$
if and only if $H$ is a subspace of $\Rn$ of dimension $n-1$, and
$\span(X\cap H)=H$.   Finally, for any facet hyperplane, let $\eta_H$
be the normal to $H$, and let $m(H)$ be the cardinality of the
vectors in $X\bks H$:
$$m(H)\eqbd m_X(H)\eqbd \#(X\bks H).$$
The three $\calI$-ideals are generated each by the polynomials
$$p_{\eta_H}^{m(H)+\eps},\quad H\in \calF(X).$$
The external ideal $\calI_+(X)$ corresponds to the choice $\eps=1$,
the central ideal $\calI(X)$ corresponds to the choice $\eps=0$,
while the internal ideal $\calI_-(X)$ corresponds to the choice
$\eps=-1$.

The Hilbert series of these three ideals are closely related
to the external activity variable of the Tutte polynomial that is
associated with $X$. We explain (and prove) this connection later.
A more rudimentary result is as follows (see Section~\ref{sec_linalg} 
for the definition of unimodularity). We denote by
$$\Pi\eqbd \C[t_1,\ldots,t_n]$$
the space of polynomials in $n$ variables, and by
$$\Pi_k^\pzero \quad (\Pi^0_k, \;\; {\rm respectively})$$
the subspace of $\Pi$ that contains all polynomials of degree $\le k$
(all homogeneous polynomials of exact degree $k$, respectively). Also, for 
any homogeneous ideal $I\subset \Pi$, we denote 
$$\ker I\eqbd \{p\in \Pi:\ q(D)p=0,\ \forall q\in I\}=
\{p\in \Pi:\ q(D)p(0)=0,\ \forall q\in I\}.  $$

Our first result provides
a very basic combinatorial connection
between the $\calI$-ideals  on the one hand and the
zonotope $Z(X)$, the integer points in it, as well as the integer points
in its interior ${\rm int} (Z(X))$.  

\begin{proposition}\label{thm_intro} Let $X \subset \R^n$ be unimodular,
and let $Z(X)$ be the associated zonotope. Then:
\item{(1)} $\dim\ker\calI_+(X)=\#(Z(X)\cap\Zn)$.
\item{(2)} $\dim\ker\calI(X)\hskip 6pt  =\vol(Z(X)).$
\item{(3)} $\dim\ker\calI_-(X)=\#({\rm int}(Z(X))\cap \Zn)$.
\end{proposition}

\nt Another related result is that
the number of (unbounded) $n$-dimensional regions in $\calH(X)$
equals $\dim \ker\calI_+(X)-\dim\ker\calI_-(X)$; this result holds
for a general $X$.
As a matter of fact, far deeper connections
between the zonotope $Z(X)$ and the $\calI(X)$-ideals are demonstrated 
in this paper: the $\calI$-ideals can be derived, each, by applying suitable 
algebro-analytic operations to a suitably chosen subset of  $Z(X)\cap \Zn$.

An important highlight of the $\calI$-ideals is that their associated
kernels can be described cleanly and explicitly  in terms of the
columns of $X$.\footnote{Unfortunately, there are no known simple 
representations for the kernels of the $\calJ$-ideals.}
Our second illustration moves in this direction and considers, for a given 
$X$, the possible use of the polynomials
$$p_Y\eqbd \prod_{y\in Y}p_y,\quad Y\subset X,$$
for the representation of the ring $\Pi$. For example,
denote 
$$\calP_+(X)\eqbd \span\{p_Y: Y\in 2^X\}.$$ 

\begin{proposition} With $\calP_+(X)$ as above,
$\Pi=\calP_+(X)\oplus \calI_+(X)$.
\end{proposition}

\nt
The above result follows directly from the fact that $\calP_+(X)$
equals $\ker \calI_+(X)$. Even more interesting decompositions are obtained
when using the $\calJ$-ideals, since these ideals are generated by
polynomials of the form $p_Y$, with $Y$ 
a (multi)subset of $X$ (or of a slightly augmented version of it).
For example, one way to express the duality between $\calI(X)$ and
$\calJ(X)$ is via the direct sum decomposition (cf.\ \S 3)
$$\Pi=\calJ(X)\oplus \ker \calI(X).$$
This decomposition corresponds to a decomposition of the power set
$2^X$: it will be shown that $\ker\calI(X)$ is spanned by 
$p_Y$, $Y\in S(X)$, with $S(X)$ a suitable subset of $2^X$. The ideal
$\calJ(X)$ is generated by the remaining polynomials $p_Y$,
$Y\in 2^X\bks S(X)$.

\medskip
Special types of zonotopal algebras are intimately connected to
group representations. The connection is particularly simple in the case
of $SL_{n+1}$-representations, since in this case the underlying
$X$ is unimodular.\footnote{Group representations are connected with
a discrete version of zonotopal algebras, that are not discussed in this 
paper. In the unimodular case, however, the discrete version coincides with 
the continuous version, which is the version studied here.} Fixing $n$, we 
let $X^k$, $k\ge 1$ be a $k$-fold multiset of the edge set\footnote{One needs
also to choose correctly the basis for $\R^n$ in the definition of the edge
set.} of a complete graph 
with $n+1$ vertices (see Example \ref{exam_edgeset}). A basic result, which 
applies to all finite-dimensional $SL_{n+1}$ representations, is that the 
character (or more precisely the Fourier coefficients of the character) of
the representation is piecewise  polynomial (see, e.g.,~\cite{Weyl_book}), with the polynomial
pieces all lying in the kernel of the ideal $\calJ(X^1)$.
Here is a rather different result.

\begin{example} Fix $n\ge 2$ and a positive integer $k$, and let
$\Gam_k$ be the irreducible $SL_{n+1}$ representation of highest weight
$(k,k,\ldots,k)$. Then there exists 
a unique polynomial $p\in\calP_+(X^k)$ whose  values on the spectrum
of $\Gam_k$ determine the character of $\Gam_k$: at each eigenvalue
$\alp$, $p(\alp)$ equals the multiplicity of the eigenvalue in $\Gam_k$.
\end{example}

This result follows directly from the theory of this paper, thanks
to the fact that the
convex hull of the spectrum of the above $\Gam_k$ is the zonotope
$Z(X^k)$. However, the connection between zonotopal algebras and group
representations extends beyond examples of this type, as the next
example makes clear. In that result, $X^k$ retains its meaning
from the previous one.   Note that, in general, the convex hull
of the spectrum of an $SL_3$-representation is not a zonotope. 

\begin{example} Let $\Gam$ be an irreducible $SL_3$ representation
of highest weight $(k+j,k-j)$, for some integers $0\le j\le k$. Then, with
$\sig$ the spectrum of $\Gam$, and for every $c\in\C^\sigma$
there exists a unique (bivariate) polynomial $p\in \calP_+(X^k)$,
such that:
\begin{itemize}
\item[(i)] $\deg p\le 3k-j$, and
\item[(ii)]  $p_{|\sigma}=c$.
\end{itemize}
\end{example}

\medskip
Let us now illustrate connections with algebraic graph theory and with
the notion of parking functions from combinatorics (see~\cite{PS_TAMS},
\cite{Yan_AAM}, \cite{KY_ArXiv}). For simplicity, we present the connection
for the edge set $X$ of
a complete graph $X$ with $n+1$ vertices; similar results
are valid for general graphical $X$. Note that here and elsewhere
$$\Z_+$$
stands for the non-negative integers (including $0$).

\begin{definition} Set $V\eqbd \{1,\ldots,n\}$. For  $r\in \Z_+^V$, and
$v\in V$, set
$V_{r,v}\eqbd \{v'\in V:\ r(v')\ge r(v)\}.$
Then $r$ is called:
\begin{itemize}
\item[(i)] An {\bf external parking function} if, for
each $v\in V$, one of  the following two conditions holds:
\begin{itemize}
\item[]  Either\ \  $\#V_{r,v}< n-r(v)+1$,
\item[]  or\ \ \ \ \ \ \ $\# V_{r,v}= n-r(v)+1$,\ \  and\ \  $r(v')=r(v)$
for $v'=\min V_{r,v}$.

\end{itemize}
\end{itemize}
\begin{itemize}
\item[(ii)] An {\bf internal parking function} if, for
each $v\in V$, one  of the following two conditions holds:
\begin{itemize}
\item[]  Either $\#V_{r,v}< n-r(v)$,
\item[]  or\ \ \ \ \ \ $\# V_{r,v}= n-r(v)$,\ \ \ \ \ \ \  and\ \  
$r(v')=r(v)$, for some $v'  \neq \max V_{r,v}.$
\end{itemize}
\end{itemize}
\end{definition} 

Parking functions define a monomial set in $\Pi$ whose monomial complement spans 
a monomial ideal. This monomial ideal ``monomizes'' a corresponding zonotopal
ideal, and the above holds for every graphical $X$; this point was already 
made explicit in~\cite{PS_TAMS} (for the central zonotopal case).
Here is a pertinent statement concerning the external case.  We use here 
$R_+(X)$ to denote the set of external parking functions of $X$.

\begin{example}  Let $X$ be the edge set of a complete graph with $n+1$
vertices. Then there exists an injection $T:R_+(X)\to 2^X$ such that
\begin{itemize}
\item The polynomials $\{p_{Tr}: r\in R_+(X)\}$ form a basis
for $\calP_+(X)$.
\item For each $r\in R_+(X)$, the monomial $t^r$ appears
(with non-zero coefficient) in the monomial expansion of $p_{Tr}$.
\end{itemize}
In particular, $\deg p_{Tr}=\sum_{v\in V}r(v)$, for every $r\in R_+(X)$.
\end{example}

Since parking functions are well known to be connected with other
combinatorial aspects of graphs, such as the number of inversions in its
spanning trees, \cite{S_1998}, results as the above draw connections with 
graph theory beyond parking functions {\it per se.\/} We study connections 
of this type in~\cite{HR_comb}.

\medskip
We now move in a completely different direction, and point out connections
between zonotopal spaces and special types of multivariate polynomial
interpolation problems. Connections of this type are at the core
of zonotopal algebras, were fully developed before for the central
case, and are well explained in the body of this paper. Here is one
illustration (cf. Section~\ref{sec_ext}).

\begin{proposition} Denote $\calZ_+(X)\eqbd Z(X)\cap \Z^n$. Then the
restriction map
$$f\mapsto f_{|\calZ_+(X)}$$
is a bijection between $\calP_+(X)$ and $\C^{\calZ_+(X)}$, provided
that $X$ is unimodular.
\end{proposition} 

Our final example is about connections of the theory
developed here with approximation theory. We recall that a (polynomial)
{\it box spline\/} $M_X$ (with $X\subset\Rn$ the given multiset) is a 
smooth piecewise polynomial function in $n$ variables supported on the 
zonotope $Z(X)$. It can be defined as the convolution product of the 
measures $M_x$, $x\in X$, with the mass of each $M_x$ uniformly distributed 
on the line connecting $0$ to $x$. One of the early key problems in box 
spline theory was to understand the properties of the polynomial space
$$\calD(X)$$
(which is defined and reviewed in Section~\ref{sec_central} here, and which is known to be)
spanned by the polynomials in the local structure of $M_X$. The ``mere''
attempt to understand the {\it dimension\/} of that space spawned  an
industry of techniques for estimating the dimension of joint kernels of 
differential and other operators (see~\cite{R_MRReview} and references 
therein). We present below a potential box spline application of our results 
that is of a different flavor. 

\begin{conjecture}\label{dynreif} Let $X$ be unimodular, and let
$$\calZ_-(X)$$
be the set of integer points in the interior of $Z(X)$. Let
$f$ be any function defined on $\calZ_-(X)$. Then, there exists
a unique polynomial $p\in \ker \calI_-(X)$ such that
$p(D)M_X$ equals $f$ on $\calZ_-(X)$.
\end{conjecture}

This conjecture follows (albeit in a somewhat non-trivial way) from
Conjecture~\ref{pminusconj}%
\footnote{%
We note that Conjecture \ref{pminusconj} implies
that point evaluation is well defined on $p(D)M_X$.}%
, but may be true even if the latter is disproved.

\subsection{\label{sec_history}Historical context}

Special zonotopal algebras (viz.\ spaces of the type $\calD(X)$
for special maps $X$) appear  implicitly in Weyl's character
formul\ae, and the connection is valid for representations
of all semi-simple  Lie algebras,~\cite{Weyl_book}. Zonotopal
spaces associated with general maps $X$ (viz.\ the spaces
$\calD(X)$) made their debut in~\cite{BH_JdAM}.  The dimension formula
for $\calD(X)$ was established in~\cite{DM_Studia} (continuous version)
and in~\cite{DM_TAMS} (discrete version). This result
was extended to non-matroidal structures by multiple authors
and in multiple ways (see~\cite{R_MRReview}). Our approach
here, in~Section~\ref{sec_central}, bypasses these developments, 
but uses in an essential way methods for bounding dimensions of such 
spaces {\it from below\/} \cite{BenR_TAMS,BR_JMAA}. The dual space 
$\calP(X)$ was introduced
independently in \cite{AS_MZ} and  in~\cite{DR_TAMS}, with the latter
containing the details concerning the construction of the homogeneous
basis for $\calP(X)$ (Section~\ref{sec_basis}). The identification
of $\calI(X)$ as the annihilating ideal of $\calP(X)$ is found
in \cite{BDR_PJM}. A chapter in~\cite{BHR} is devoted to the study of
these and other related algebraic aspects of box spline theory.
Newer treatments of the central case are presented in~\cite{DP_ArXiv} 
and the book~\cite{DP_book},  where several aspects of the central 
algebra are re-explored and extended, including its relations with 
modules over the Weyl algebra and $D$-modules, as well as with
toric arrangements and their cohomology (see also~\cite{DP_note}). 
A nice connection between 
the space $\calD(X)$ and the cohomology of toric hyperK\"ahler 
varieties is described in~\cite{HS_DM} via the so-called Volume Polynomials 
shown to span $\calD(X)$ as a $D$-module (for subsequent developments,
see~\cite{HP_T}, \cite{HS_PAMS}). 

Our interest in extending the theory of zonotopal spaces beyond the
central pair was stimulated by discussions we had in the mid 1990's
with N. Dyn and U. Reif, concerning the possibility of a result  
{\it \'a la\/} Conjecture~\ref{dynreif},   and was enhanced  by 
discussions we had a few years later with F. Sottile, who
pointed out to us connections between our external theory and the work
of~\cite{PSS}. Our delay in publishing this theory was primarily  
caused by inherent difficulties we encountered in the internal study 
due to the absence of a ``canonical'' basis for $\ker\calI_-(X)$. 
We believe  that the theory as presented here alleviates ramifications 
of this hurdle to the extent possible.

As we alluded to above,
the novelty of this paper lies exclusively in the theory of the internal and 
external algebra, whose foundations we develop here, as well as in pointing
out various connections of this theory with other fields -- most
notably, enumerative combinatorics and representation theory  (see
Sections~\ref{sec_ext} and~\ref{sec_int}, also Section~\ref{sec_general}). 
It should be mentioned that our second task is by no means completed 
in this paper, due to the multitude and richness of those connections.
A description of combinatorial connections alone is a subject of our
forthcoming paper~\cite{HR_comb}, currently in preparation. 

We hope that this paper will offer a new perspective and new tools 
to researchers working in algebra, analysis and combinatorics, along 
with a glimpse into exciting developments yet to come.

\noindent
{\bf Added February 2011.} We take this opportunity to comment on various
developments relevant to the subject of this paper that have taken place 
since August 2007, when the preprint of this paper was initially posted on the arXiv. 

The most significant of those is the completion and publication
of a monograph~\cite{DP_book} by C.~de Concini and C. Procesi, which 
gives a detailed account of  the algebraic connections between 
polytope geometry and combinatorics, differential and difference 
equations, and box spline theory. 
Another significant achievement is a series of  
papers~\cite{dCPV1,dCPV2,dCPV3} of C. de Concini, C. Procesi 
and M. Vergne that trace connections between
vector partition functions, zonotopal spaces, and
indices of transversally elliptic operators.
This body of work provides a novel framework to understand
the index map of such operators, the setting studied by
M. F. Atiyah and I. M. Singer.

In another exciting work \cite{StuXu08} partly motivated by this paper,
B. Sturmfels and Z. Xu study sagbi bases of Cox-Nagata rings and
generalize the graded zonotopal algebras to the so-called zonotopal Cox rings.
F. Ardila and A. Postnikov~\cite{AP08} investigate the setting of power 
ideals, which generalize the zonotopal polynomial
ideals from this paper.
Together with Z. Xu, we construct a whole hierarchy
of zonotopal spaces interpolating between the internal and external 
spaces in~\cite{HRX09}. M. Lenz~\cite{Lenz10} further generalizes 
a number of results from~\cite{HRX09} and~\cite{AP08}.

A brief note is in order regarding the status of 
Conjecture~\ref{pminusconj}. F. Ardila and A. Postnikov claim a 
proof of this conjecture in~\cite{AP08}. However, M. Lenz points
out~\cite{Lenz10} that the proof contains a mistake.
F. Ardila and A. Postnikov have recently  
corrected their proof, according to~\cite{Ardila_private}.

A. Berget \cite{Berget09} clarifies connections between the zonotopal
spaces, the fundamental work of P.~Orlik and H. Terao on hyperplane 
arrangements \cite{OT94,Terao02}, and the algebras studied by D. Wagner
 in~\cite{Wagner99}. 
L. Moci \cite{Moci09} generalizes our results on hyperplane
arrangements to toric arrangements (see also~\cite{DAM11}).  
One of Moci's most interesting results shows that, in his 
(discrete) setting, the connection between the dimension of 
the space and the volume of the zonotope remains valid even 
if the matrix $X$ is not unimodular.

Finally, several combinatorial connections of this work to the Ehrhart 
polynomial, Tutte polynomial, and other counting functions on graphs
and matroids are by now worked out \cite{AP08,Berget09,DAM11,Des10,HR_comb}.
The most important combinatorial connection to the Tutte polynomial
is clarified in~\cite{AP08} and \cite{Berget09} (for toric generalizations,
see~\cite{Moci09,DAM11}). It deserves to be mentioned explicitly here as well.

A comparison of  the grading of the spaces  $\calP(X)$, $\calP_-(X)$ and
 $\calP_+(X)$ (or, equivalently,  their duals) described in 
Sections~\ref{sec_basis},  \ref{sec_basis+}, \ref{sec_basis-}
with the definition of the Tutte polynomial $T_X$ 
of the matroid of $X$ via the external and internal activity 
demonstrates that the Hilbert series ${\rm Hilb}_X$, ${\rm Hilb}_{X,-}$ and 
${\rm Hilb}_{X,+}$  of the spaces $\calP(X)$, $\calP_-(X)$ and $\calP_+(X)$,
respectively, 
are the following  specializations of the Tutte polynomial
(cf.~\cite[Corollary 4.13, Proposition 4.14, Proposition 4.15]{AP08}):
\begin{proposition}\label{prop_tutte}
Let $X \subset \R^n$ be of full row rank $n$ and consist of $N=\# X$
vectors. Let $T_X$
denote the Tutte polynomial of the matroid of $X$ and let
${\rm Hilb}_X$, ${\rm Hilb}_{X,-}$ and ${\rm Hilb}_{X,+}$  denote the Hilbert
series of the spaces   $\calP(X)$, $\calP_-(X)$ and $\calP_+(X)$,
respectively. Then   
\begin{eqnarray*} {\rm Hilb}_X(q)=q^{N-n}T_X\left(1,{1\over q}\right), \quad
{\rm Hilb}_{X,+}(q)=q^{N-n} T_X\left(1+q,{1\over q}\right), \quad
{\rm Hilb}_{X,-}(q)=q^{N-n} T_X\left(0,{1\over q}\right). 
\end{eqnarray*}
\end{proposition}
As the main result of~\cite{Berget09}
shows, we can in fact also obtain the full bivariate 
Tutte polynomial as the Hilbert series of the space 
${\calP}_+(X)$ if we equip that space
with  a suitable bi-graded structure.

\subsection{\label{sec_layout}Layout of this article}

The paper is organized as follows. Section~\ref{sec_prelim} contains
background results that will be used in the rest of the paper. This
section is subdivided into five subsections: Section~\ref{sec_linalg}
is devoted to  linear algebra and matroid theory, Section~\ref{sec_hyparr}
to hyperplane arrangements, Section~\ref{sec_zono} to zonotopes, 
Section~\ref{sec_least} to polynomial interpolation and 
Section~\ref{sec_ideals} to polynomial ideals and their kernels.

The bulk of the paper is in Sections~\ref{sec_central}, \ref{sec_ext}
and \ref{sec_int}. Those three sections are made  parallel to each 
other, with two subsections in each, the first containing  main theory
and the second discussing grading, the Hilbert series, and homogeneous
bases for the polynomial spaces in question. While the material in
Section~\ref{sec_central} is known, we feel it is crucial to present it
in this way here, for the rest of the paper to be much more easily 
understandable, as well as for the streamlined approach itself. 
The paper ends with Section~\ref{sec_conclude} containing a few additional
remarks and conjectures.

\section{\label{sec_prelim}Preliminary results}

\subsection{\label{sec_linalg}Linear algebra}

Consider a finite multiset $X \subset \R^n \backslash \{0\}$ of full 
rank $n$ and of size $\# X$. At times, we will associate $X$ with some 
full ordering.  In this case, we may consider the vectors in $X$ to 
comprise the columns of an $n$-by-$\# X$ matrix, which we will still 
denote by $X$. Such a multiset (or a matrix) $X$ gives rise to a {\it 
linear matroid\/} (see, e.g.,~\cite{Oxley}) via the standard convention 
that  the independent sets  of the matroid are exactly the linearly 
independent  subsets of the columns of $X$.

We now single out three sub-collections of the power set $2^X$ that will
play a  crucial role in this paper. The reader may notice right
away that all three are defined in  purely matroidal terms.
The first of the three is the set $\B(X)$ of all bases of $X$:
$$\B \eqbd \B(X) \eqbd  \Set{ B\subset X \,:\, B \mbox{ is a 
basis for } \R^n}.  $$ 
The second is the collection $\I(X)$ of all independent subsets of $X$:
$$ \I \eqbd  \I(X) \eqbd \Set{ I \subset X ,\:\, 
 I \mbox{ is independent in } \R^n}.$$
Note that the empty set is independent. The third is the set of internal 
bases, and is defined in the sequel.

Clearly, $\B(X)\subset \I(X)$, and the inclusion is proper.
Nonetheless, it is convenient to consider the independent sets 
as full-rank bases, too.  To this end, we choose 
a fixed basis $B_0$ of $\R^n$,  and append $B_0$ to $X$:
$$X'\eqbd  X \cup B_0.$$
We then impose some arbitrary, but fixed, ordering $\prec$ on $B_0$, and 
associate  each $I \in \I(X)$ with $\ex(I)\in \B(X')$ which is the greedy
completion of $I$ to a basis, using the elements of $B_0$, i.e.,
$b\in \ex(I)$ if and only if $b\in I$ or else $b\in B_0$ and
$$b\not\in \spam\{I\cup\{b'\in B_0: b'\prec b\}\}.$$
That creates a 1-1 map from $\I(X)$ into $\B(X')$. The range of this
extension map is denoted 
by $\B_+(X)$:
$$\B_+(X) \eqbd  \Set{\ex(I)\in \B(X') \,:\, 
  \ I \in \I(X)}.$$
We refer to the bases in $\B_+(X)$ as the {\it external bases\/} of
$X$. Note that every basis of $X$ is external directly from the definition, 
but not every external basis of $X$ is a basis of $X$.

Next, we define the notion of an internal basis.  To this end, we assume 
to be given an order $\prec$ on $X$. A vector $b\in B$ in a basis $B \in 
\B(X)$ is said to be {\it internally active\/} in $B$ if $b$ is the last element 
in $X\backslash H$, where $H \eqbd  \spam (B\backslash b)$:
\begin{equation}\label{defintact}
b\succ x,\ \forall x\in X\bks (H\cup b).
\end{equation}
A basis
$B$ that contains {\it no\/} internally active vectors is said to be
an {\it internal basis\/}. We denote 
$$\B_-(X) \eqbd  \Set{ B\in \B(X) \,:\,   B \mbox{ is an internal basis}}.$$
It is obvious that the notion of an internal basis depends on the ordering.
In fact, assuming that the last $n$ vectors of $X$ form a basis $B_1$,
only the ordering within $B_1$ counts here, since, whatever $B\in\B(X)$
we choose, only the vectors in $B\cap B_1$ can be internally active in $B$.
Thus
$$\B(X\bks B_1)\subset \B_-(X)\subset\B(X)\subset \B_+(X)\subset\B(X\cup
B_0).$$
We will see later  that the number of internal bases is
independent of the ordering of $X$.

We say that $X$ is {\it unimodular\/} \label{p_unimod} if $X \subset \Z^n$ and
$$\forall B\in \B(X),\ \mbox{span}_\Z B = \Z^n 
 \ (\iff \mbox{det} (B) = \pm 1).$$


\begin{example}\label{exam_edgeset}{\bf 
[the edge set of a graph].  \/  \rm Let $G$ be a connected undirected graph 
with $n+1$ vertices $V=\{v_i\}_{i=0}^n$. Let $e_0\eqbd 0$. Let $(e_i)_{i=1}^n$
 be a basis for $\R^n$. Identify an edge $e_{ij}$ that connects the vertices 
$v_i$ and $v_j$ ($i<j$) with the  vector $e_i-e_j\in\R^n$. With this identification, 
one chooses $X$ to be the edge set of $G$. Note that the edge (multi)set 
$X$ of a graph is always unimodular (assuming, say, that $(e_i)$ is the
standard basis. Otherwise, ``unimodularity'' here is with respect to the
lattice spanned by the basis.) The corresponding matroid is called 
{\it graphical.\/} A particular interesting example is when $G$ is chosen to
be a {\it complete graph,\/} i.e., a graph in which every pair of vertices is
connected by exactly one edge.}
\end{example}

\nt {\bf Remark.\/}
Although it is not obvious, there is a certain level of symmetry in the
definition of external bases and internal ones. To demonstrate this point,
let us assume that $X$ is graphical, and let $B_1\eqbd (e_i)_{i=1}^n$. Assuming
that $B_1\in\B(X)$ (which means that there is an edge between $v_0$ and each
of the other vertices), we place $B_1$ last in $X$, and order its vectors
according to the enumeration of the vertices ($e_i\prec e_j$
iff $i<j$). Using this order to define $\B_-(X)$, one finds that 
$B\in \B(X)$ is internal if $B\in \B(X\bks B_1)$. Otherwise, $B\bks B_1$
defines a partition $(V_0,\ldots,V_k)$, $v_0\in V_0$,
on the vertex set $V$, with the $k$
vectors in  $B\cap B_1$ connecting $v_0$ to each of $V_1,\ldots,V_k$.
The basis $B$ is then internal if and only if, for $i=1,\ldots,k$, the edge
 in $B\cap B_1$ that connects $v_0$ and $V_i$ is {\it not\/} connected to the 
{\it maximal\/} vertex of $V_i$.

Now, let us append another copy of $B_1$ to $X$. We call this
copy $B_0$, and denote $X' \eqbd X\cup B_0$. We
 retain the
order on $B_0$ as above and use this external copy $B_0$
to define $\B_+(X)$. We then need  to determine what the greedy extension
$\ex(I)$ of $I\in\I(X)$ should be. Again, each such $I$ determines a partition
$(V_0,\ldots,V_k)$ as before. The greedy extension is performed by connecting,
for $i=1,\ldots,k$, the vertex $v_0$ to the {\it minimal\/} vertex in $V_i$.
\eop

\subsection{\label{sec_hyparr}Hyperplane arrangements}

Recall from the introduction the definition 
$$ \Pi \eqbd \C[t_1, \cdots, t_n],$$   
as well as the notations $\Pi_k^\pzero$ and $\Pi_k^0$.
We first associate each direction $x\in X$ with a constant
$\lam_x\in \R$, and define an affine polynomial $p_{x,\lambda} \in \Pi$: 
$$p_{x,\lambda} : \R^n\to \R : t \mapsto x\cdot t - \lam_x.  $$
The $X$-{\it hyperplane arrangement\/} $\calH(X,\lam)$ is determined
by the zero sets of the above polynomials, viz., by
$$H_{x,\lambda} \eqbd  \Set{t \in \R^n \,:\, p_{x,\lambda}(t) =0},\quad x\in X.$$
We will usually assume that $\lam$ is chosen {\it generically,\/} 
i.e., so that 
the intersection of any collection of $n+1$ hyperplanes  is empty.  
Note that different choices of $\lam$ may result in 
hyperplane arrangements with different geometries. Of particular
interest are the following three geometric characteristics of the 
hyperplane arrangement:
\bd
\item[1.] $V(X,\lam)$: the set of vertices
\item[2.] $CC(X,\lam)$: the set of $n$-dimensional connected components
\item[3.] $BCC(X,\lam)$: the set of $n$-dimensional bounded connected components
\ed

As a reader of this article should observe, the set $V(X,\lam)$
is analyzed in Section~\ref{sec_central}; however, the set $CC(X,\lam)$ appears nowhere in this
paper past the current location. The reason is mainly technical:
the tools that we introduce and employ allow us to study zero-dimensional
sets. We bypass this limitation by associating $CC(X,\lam)$ and $BCC(X,\lam)$
with suitable supersets and/or subsets of $V(X,\lam)$, and utilize to this end
the notions of external and internal bases.  It is thus worth mentioning
the following known facts. 

\begin{result}   \label{thm_count}
For any generic hyperplane arrangement determined by a multiset $X$, 
\begin{eqnarray*}  \# V(X,\lam) & = & \# \B(X), \\
 \# CC(X,\lam) & = & \# \B_+(X) \; ( = \; \#\I(X) ), \\ 
 \# BCC(X,\lam) & = & \# \B_-(X). \end{eqnarray*}
\end{result}

The result shows in addition that the number of objects of each type is a 
geometric invariant of generic arrangements.
In this connection, it is worthwhile to point out the relevance of the
(univariate) {\it Ehrhart polynomial\/}\footnote{Strictly speaking, this is 
the true Ehrhart polynomial of $X$ only when $X$ is unimodular (see, 
e.g.,~\cite{S_DIMACS}).}:
$$ E_X(t) \eqbd  \sum_{I\in\I(X)} t^{\#I}.$$

It is known that~(see~\cite{BeckR_2007})
$E_X(1) =\#\I(X)=\# \B_+(X) $, and that $E_X(-1) =(-1)^n \# \B_-(X) $.
The first equality is a triviality; the second
one can be proved by induction on $n$. It shows that
$\#\B_-(X)$, indeed, is independent of the order on $X$.

\subsection{\label{sec_zono}Zonotopes}

Let us now consider $X$ as a map:
$$X: \R^X \to \R^n : t \mapsto \sum_{x\in X} t_x x.$$
Then the {\it zonotope of $X$\/} is defined as the image of the unit cube
under this map $$ Z(X) \eqbd  X ([0,1]^X).$$
Assuming $X$ to be unimodular, we have the following formul{\ae}  
for the volume of $Z(X)$, the number of integer points in $Z(X)$,
and the number of integer points in the interior of $Z(X)$, respectively:
\bd
\item[1.] $\mbox{vol}(Z(X)) = \#\B(X)$,
\item[2.] $\# (Z(X)\cap \Z^n) = \# \B_+(X)= \#\I(X), $
\item[3.] $\# (\mbox{int}(Z(X))\cap \Z^n) = \#\B_-(X)$.
\ed
Every zonotope $Z(X)$ is a disjoint (up to a nullset) union 
of the translated parallelepipeds~\cite{S_DIMACS,BHR}:
$$t_B+Z(B),\quad B\in\B(X).$$
The translation $t_B\in \R^n$ equals $\sum_{x\in X(B)}x$, with
$X(B)$ a suitable subset of $X\bks B$.  There are multiple ways of 
choosing these translations, hence there are multiple tilings of 
the zonotope.  
A canonical approach to obtaining a tiling is based on ordering
$X$ (in any way). Each such ordering corresponds to a different geometry
on the hyperplane arrangement. In this duality, the vertices
of the hyperplane arrangements are associated with the parallelepipeds
that tile the zonotope, the bounded regions of the arrangement correspond
to the vertices of the parallelepipeds that are interior to the zonotope, 
while the unbounded regions of the arrangement correspond to the vertices 
on the boundary of the zonotope. Thus, for example, the number of vertices 
of a connected region of the arrangement must agree with the 
number of parallelepipeds that intersect at the corresponding
``lattice point'' of the zonotope.  This geometric duality is well known and 
is discussed, e.g., in~\cite{S_notes, S_1998, BDDPS_2005, BeckR_2007}. 
A reader who is interested in the above-mentioned geometric duality
may wish to revisit the discussion here after reviewing the construction of a 
homogeneous basis for  $\calP(X)$ in Section~\ref{sec_basis}.

\subsection{\label{sec_least}The least map and polynomial interpolation}

Given a power series $f$ in $n$ variables
$$ f = f_0 + f_1 +  f_2 +  \cdots, $$
where $f_j$ is a homogeneous polynomial of degree $j$,
and  define the {\it least map\/} via $f\mapsto f_\downarrow$ by
$$f_\downarrow \eqbd f_j,\quad f_j\not=0, \;\; f_i=0,\;\; \forall i<j.$$  
In other words, $f_\downarrow$ is the first non-zero term in the above 
expansion of $f$. We adopt the convention that  $0_\downarrow\eqbd  0$. 
For a collection $F$ of functions analytic at the origin, we define 
$$  F_\downarrow \eqbd \spa \{ f_\downarrow : f \in F \}. $$ 


The least map plays an important role in polynomial interpolation, as 
was shown in \cite{BDR_PJM, BR_CA, BR_JMAA}.  Here are the 
details.  A pointset $\sigma\subset\C^n$ is called {\it correct\/} for a 
polynomial space $P\subset \Pi$ if the restriction map $p\mapsto p|_\sigma$ is 
invertible (as a map from $P$ to $\C^\sigma$) i.e., 
if  interpolation from $P$ at the points of $\sigma$ is {\it correct\/}
(the latter means that the interpolating polynomial exists and is
unique).

With $\sig$ a finite subset of $\R^n$, consider
the  point evaluation functional
$$\del_\alp : \Pi \to \C \,:\, p \mapsto p(\alp),$$
and define
$$\Lam \eqbd  \span \Set{\del_\alp \,:\, \alp \in \sig}.$$
Then $\Lam$ is  a subspace of the dual space $\Pi'$ of $\Pi$.   
Given $\sig$, the correctness of a space $F\subset \Pi$ for interpolation
on $\sig$ is equivalent to the isomorphism
$$ \Lam |_F  \cong F',$$
where $F'$ is the dual space of $F$. 

Now, associate $p\in \Pi$ with a differential
operator: $$ p(D)\eqbd p \gp{ \frac{\partial}{\partial t_1}, \cdots,  
 \frac{\partial}{\partial t_n}}.$$
(In particular, if $p(t) = x\cdot t$, \/ $x \in \R^n\backs \{0\}$, then 
$p(D)$ is the directional derivative in the $x$-direction.)
Then,
given  a polynomial $p$ and a formal power series $f$, define their pairing
$\inpro{p,f}$ as 
\begin{equation} \inpro{p, f} \eqbd  (p(D) f)(0). \label{inner_prod} 
\end{equation}
The functional $\del_\alp$ is represented using  this pairing 
by the exponential
$$e_\alp : \R^n \to \R \,:\, t \mapsto e^{\alp \cdot t},\quad {i.e.,}  \qquad 
\inpro{p, e_\alp} = p(\alp) = \del_\alp p.$$
Thus the space $\Lam$ is represented by the exponential space
$$\Exp (\sig) \eqbd  \span \Set{ e_\alp \,:\, \alp \in\sig}.$$
Finally, we define 
$$\Pi(\sig) \eqbd  \span \Set{ f_\downarrow \,:\,  f\in \Exp(\sig)}.$$

We now check that the dimension of $\Pi(\sig)$ is exactly $\# \sigma$. 
Let $T_j:\Exp(\sig)$ be the $j$th degree Taylor expansion; i.e.,
for $f\in \Exp(\sig)$, $T_j f $ is the $j$-th degree Taylor
expansion of $f$ at $0$.
Note that $\deg(f_\downarrow)=j$ if and only if $f\in\ker T_{j-1}\backslash 
\ker T_j$.  Thus, 
with $T'_j$ the restriction of $T_j$ to $\ker T_{j-1}$,
$$\dim {\rm span}\{f\least: f\in \Exp(\sig),\ \deg(f\least)=j\}
=\rank T'_j=\dim\ker T_{j-1}-\dim\ker T_j.$$
Summing from $j=0$ to $\infty$ (where $T_{-1}\eqbd 0$), we obtain 
$$ \sum_{j=0}^\infty \Big( \dim(\ker T_{j-1}) -\dim(\ker T_j)\Big)
=\dim(\ker T_{-1})=
\dim \Exp(\sig)=\#\sig.$$
Here we used the fact that every finite set of exponentials is linearly
independent.

\noindent
Now, for any analytic function $f\not=0$, 
$$ \inpro{f_\downarrow, f} \neq 0.$$
This means that there exists no $f\in \Exp(\sig)\backs \{0\}$ that satisfies
$$\inpro{p,f} =0, \quad \forall p\in \Pi(\sig).$$
Thus, the spaces $\Exp(\sig)$ and $\Pi(\sig)$ serve as duals
of each other.  In summary:

\begin{result}[\cite{BR_CA}] \label{res_pi(sig)}
The map $\Exp(\sig)\to \Pi(\sig)'$ defined by $f \mapsto \langle \cdot,f   
\rangle$  is an isomorphism, and  the set 
$\sig$ is correct for the space $\Pi(\sig)$. In particular,
$\dim \Pi(\sig)=\#\sig$. 
\end{result}


Now, given any non-zero polynomial $p\in \Pi$, let $p_\uparrow$ be the highest 
degree homogeneous polynomial in $p$, i.e., $p_\uparrow$ is
homogeneous, and $\deg(p-p_\uparrow)<\deg p$.  Likewise, define, for
$F\subset \Pi$,
$$F_\uparrow \eqbd  \span \Set{ f_\uparrow \,:\, f \in F}. $$
This defines the so-called {\it most map.\/} The following result
describes the interaction of the least map $(\cdot)_\downarrow$
and the most map $(\cdot)_\uparrow$:

\begin{result}[\cite{BR_CA}]
Let $p\in \Pi$ be fixed. Then
\begin{equation}
 (p(D) \gp{ \Exp(\sig)})_\downarrow \supset p_\uparrow(D) (\Pi(\sig)).
\label{exp_least}
\end{equation}
\end{result}

\nt {\bf Proof. \/} 
Let $p$ be a polynomial and let $f$ be an analytic function.  
Set $p \bdeq p_\uparrow + q$, then $\deg q < \deg p$. 
For $f\in \Exp(\sig)$, set $f \bdeq f_\downarrow + g$.
Then we have 
$$p(D)f  = p_\uparrow(D) f_\downarrow + q(D) f_\downarrow + 
p_\uparrow(D) g + q(D) g.  $$
Note that $p_\uparrow(D) f_\downarrow$ is the lowest order term
in the right hand side. Hence $p_\uparrow(D) f_\downarrow = 
(p(D)f)_\downarrow$
unless $p_\uparrow(D) f_\downarrow=0$.
Applying this observation to an arbitrary function $f\in \Exp(\sig)$,
we conclude that~(\ref{exp_least}) holds.   \qed

This theorem can be used as follows: suppose that we have a
homogeneous polynomial $r$, and we would like to understand the action
of $r(D)$ on $\Pi(\sig)$. Then it makes sense to find an inhomogeneous 
polynomial $p$ such that (i) $p_\uparrow=r$, and (ii) $p$ vanishes
at as many points of $\sig$ as possible. Indeed, one easily verifies that
$$p(D)\Exp(\sig)=\Exp(\sig\backslash Z_p),$$
with $Z_p$ the zero-set of $p$. In particular, $p(D)$
annihilates $\Exp(\sig)$ iff $p$ vanishes on $\sig$. In that latter case,
we obtain the following corollary from the previous result:

\begin{corollary}[\cite{BR_CA}]  \label{cor_annihil}
If $p$ vanishes on $\sig$, then $p_\uparrow$ annihilates $\Pi(\sig)$.
\end{corollary}

\begin{corollary}[\cite{BR_CA}]  \label{cor_duala}
Let $\sig\subset\R^n$ be finite, and let $P$ be a homogeneous
subspace of $\Pi$.  Then $\sig$ is correct for $P$ if the map
$$p\mapsto \inpro{\cdot,p}$$
is a bijection between $P$ and $\Pi(\sig)'$.\footnote{%
The converse of this result is true, too, once we assume in addition
that $\dim( P\cap\Pi_j)\ge \dim(\Pi(\sig)\cap\Pi_j)$, for all $j$.}
\end{corollary}

\nt {\bf Proof. \/}  We may assume that $\dim P=\#\sig$, since
otherwise the bijection cannot hold.  Now, if $\dim P=\#\sig$, and
$\sig$ is not correct for $P$, then some $p\in P$ vanishes on 
$\sig$, hence $p_\uparrow(D)$ annihilates $\Pi(\sig)$, and,
in particular, $p\most\perp\Pi(\sig)$. Since
$P$ is homogeneous, $p_\uparrow\in P$, hence the bijection does not hold.
\eop

\nt These results can be used to prove dimension formul{\ae} for 
polynomial spaces of interest, which we now introduce. 
Given a vector $\lambda$ indexed by $X$, recall that we 
associate  each  $x\in X$ with an affine  polynomial 
$$p_{x,\lambda} : \R^n\to \R : t \mapsto x\cdot t - \lam_x.$$
For simplicity, we denote the linear polynomial $p_{x,0}$ by $p_x$. 
For a multi-subset $Y\subset X$, define
\begin{equation}
 p_{Y,\lambda} \eqbd  \prod_{y\in Y} p_{y,\lambda},  \qquad
p_{Y} \eqbd p_{Y,0}. \label{def:p_Y}
\end{equation}
Let $V(X,\lam)$ denote the vertex set of the corresponding 
$X$-hyperplane arrangement $\calH(X,\lam)$.
Recall that $\# V(X,\lam)=\# \B(X)$ for a generic $\lambda$.
Moreover, if $\lambda$ is generic,  then there is a natural bijection 
$B\mapsto v_B$ between $\B(X)$ and $V(X,\lam)$,
where each $B\in \B(X)$ is mapped to the unique common zero $v_B$
of $\{ p_{y,\lam}:\ y\in B \}$.  
This implies that each subset $\B'$ of $\B(X)$ is associated 
in a unique way with $V' \subset V(X,\lam)$. We define
$$\calD_{\B'} (X) \eqbd \Set{ f\in \Pi \,:\, 
p_Y(D)f=0,\quad \forall Y\subset X \;\; \mbox{s.t.}\;\; Y\cap B \neq
\emptyset, \ \forall B \in \B'}.$$ 
Then we have the following results, \cite{BR_JMAA}:
\begin{theorem}\label{thm_db'} Let $X \subset \R^n \backslash \{0\}$ be a 
finite multiset of full rank $n$. Then, for any subset $\B'$ of\/ $\B(X)$, 
$$ \Pi(V') \subseteq \calD_{\B'} (X), $$
with $V'$  the vertices that correspond to $\B'$
in any generic $X$-hyperplane arrangement
$\calH(X,\lam)$.
\end{theorem}
 \begin{corollary}\label{cor_db'} In the setting of Theorem~\ref{thm_db'},
$$ \dim \calD_{\B'} (X) \ge \# \B'. $$
\end{corollary}

\nt{\bf Proof of Theorem~\ref{thm_db'} and Corollary~\ref{cor_db'}.\/}
Let $\lambda$ be generic. Note that, for $Y\subset X$ and $B\in\B(X)$, 
$p_{Y,\lambda}(v_B)=0$ if and only if
$Y\cap B\not=\emptyset$. Now, set $V'\eqbd \{v_B:\ B\in\B'\}$, and
let $Y$ be a multi-subset of $X$ such that
$$Y\cap B \neq \emptyset \quad \forall B \in \B'. $$
Then we conclude that $p_{Y,\lambda}$ vanishes on $V'$.
Then, by Corollary~\ref{cor_annihil}, we have 
$$ (p_{Y,\lambda})_\uparrow(D)(\Pi(V')) = 0.$$
Since $(p_{Y,\lambda})_\uparrow=p_Y$,
we conclude that
$$\Pi(V') \subseteq \calD_{\B'} (X).  $$
But, 
$$ \dim \Pi(V') = \# V' =\# \B'. $$
\eop

We will apply these results three times in this paper: once with
respect to $\B'\eqbd \B(X)$, then with respect to $\B'\eqbd \B_+(X)$,
and finally with respect to $\B'\eqbd \B_-(X)$. In all these cases,
we will show that equality holds in  Corollary~\ref{cor_db'},
hence that $\Pi(V')=\calD_{\B'}(X)$.

\subsection{\label{sec_ideals}Polynomial ideals and their kernels}

Here we state, for the convenience of the reader, some basic results from 
commutative algebra, which will be used in the rest of the paper. The 
majority of these results  are standard and can be found 
in commutative algebra textbooks, e.g.,~\cite{Eisenbud} or~\cite{Kaplansky}. 
The fact that kernels of polynomial ideals can be synthesized  from finitely 
many localizations can be found  in~\cite{Lefranc_CRAS};  while this result
is non-trivial, we will need only the result for the simpler special
case of zero-dimensional ideals. Some of the actual 
presentation here follows~\cite{BR_JMAA}.

Let $I$ be an ideal  in the ring $\Pi \eqbd  \C[t_1, \cdots, t_n].$
If $I$ is generated by a set $L\subset \Pi$, this will be denoted as
$I = \Ideal (L)$. 
The {\it codimension\/} of $I$ is the dimension of the quotient space
$\Pi/I$ or, equivalently, the dimension of its {\it annihilator\/} 
$$ \{  \mu \in \Pi'\; : \;  \mu f =0 \quad \forall f \in I   \}.  $$
Since the dual space $\Pi'$ can be realized as the space
$\C[[t_1,\ldots,t_n]]$ of formal power series, the codimension 
of $I$ is also the dimension of the orthogonal complement of $I$ in
$\C[[t_1,\ldots,t_n]]$ with respect to the pairing~(\ref{inner_prod}),
i.e., $  \langle f , g \rangle \eqbd  (f(D) g)(0). $
The {\it variety\/} $\var(I)$ of $I$ 
$$ \var(I) \eqbd \{ \theta \in \C^n \; : \;  p(\theta)=0 \quad \forall p\in I\}$$
forms a subset of the annihilator of $I$, where evaluation at a 
point $\theta\in \var(I)$ is realized by the exponential $e_\theta$.
The {\it kernel\/} of $I$ is defined as
$$ \ker I \eqbd  \span
\Set{ e_{\alp} p \, : \, \alp \in \var(I), \;\; p\in \Pi \;\; 
\mbox{ s.t. } \;\; \inpro{e_{\alp} p, q} =0, \quad \forall\, q \in I}. $$
With that definition, the kernel $\ker I$ is
{\it total\/} in the sense that
$$(\ker I)\perp \eqbd  \Set{q\in \Pi \,:\,  
\inpro{k, q} =0, \quad \forall\, k\in \ker I} =I.$$
 
An ideal $I$ is called {\it zero-dimensional\/} if its variety is finite.
In that case, each of the multiplicity spaces 
$$(\ker I)_\alp\eqbd \{p\in \Pi: e_\alp p\in \ker I\}$$
is finite-dimensional, hence $\ker I$ is a finite-dimensional space
of exponential polynomials:
$$\ker I \eqbd  \sum_{\alp\in \var(I)}(\ker I)_\alp.$$
Furthermore, we have then that
$$\dim (\ker I) = \dim \Pi \slash I,$$
hence, in particular, $I$ is of finite codimension.
Defining
$$(\ker I)_\downarrow\eqbd \span\{f_\downarrow: f\in \ker I\},$$
we have
\begin{result}
If  $I$ is a zero-dimensional ideal, then 
$$\Pi = I \oplus (\ker I)\least. $$
\end{result}



The homogenization  of the kernel of an ideal via the least map 
$(\cdot)_\downarrow$ is dual to the homogenization of the ideal itself 
via the most map $(\cdot)_\uparrow$.  Here are the details. 
Given an ideal $I$,  we define 
$$ I\most \eqbd  \spam \Set{p\most \,:\, p \in I}. $$
$I\most$ is a {\it homogeneous\/} ideal, i.e., is generated
by homogeneous polynomials. 
We have:

\begin{result}
For a zero-dimensional ideal $I$, the following properties hold:
\bd
\item[1.] $\var(I\most) = \Set{0}.$
\item[2.] $\ker (I\most)$ is a finite-dimensional polynomial space.
\item[3.] $\ker (I\most)=\{p\in \Pi: q(D)p=0, \ q\in L\}$, with
$L\subset\Pi$ any set that generates $I\most$.
\item[4.] $\ker(I\most) = (\ker I)\least.$
\item[5.] $ (I\most) \oplus \ker(I \most) = \Pi.$
\ed
\end{result}

Finally, if $I$ is $0$-dimensional, and $F$ is a polynomial space,
then the relation $F+I=\Pi$ implies that $\dim F\ge \dim\Pi/I$ with
equality iff the above sum is direct. Hence we have:

\begin{result} \label{res_sum}
Let $I$ be a zero-dimensional ideal,
and let $F$ be a linear subspace of $\Pi$. If 
\begin{equation} F + I = \Pi,  \label{F+I=Pi}
\end{equation}
then $\dim F \geq \dim \ker I$. Moreover, if\/ $\dim F = \dim \ker I$,
then the sum~(\ref{F+I=Pi}) is direct.  
\end{result}

\section{\label{sec_central}Central zonotopal algebra}

\subsection{\label{main_central}Main results}

We have mentioned the geometric duality between the hyperplane arrangement
and the zonotope that are associated with the multiset $X$. The focus
of this paper is on an algebraic counterpart of that duality. We 
discuss in the paper three pairs of finite-dimensional polynomial spaces, all
of which can be alternatively described as kernels of certain 
zero-dimensional ideals. Each polynomial space  will be shown to be 
a dual space of its pair-mate via the map $p\mapsto \langle p, \cdot \rangle$
where $\langle \cdot, \cdot \rangle$ is our pairing~(\ref{inner_prod}). 

The first pair will be referred to as the {\it central pair\/} of $X$.
The space $\calP(X)$ below is the central space of the zonotope
$Z(X)$, while the space $\calD(X)$ is the central space of the
hyperplane arrangement $\calH(X,0)$. The theory of this pair of polynomial 
spaces was developed in the 80s and 90s in~\cite{AS_MZ, AS_Izv, AS_Uspekhi, 
BenR_TAMS, BR_CA, BR_JMAA, BR_MC, BR_MZ, BDR_PJM, DM_LAA, DM_Studia, 
DM_Adv, DM_MN, DM_AAM,  DR_TAMS, DR_AA}. We have discussed some historical
aspects of this theory in the Introduction, and will discuss the history 
of specific results in more detail towards the end of this section.
The section contains a streamlined and abbreviated theory of
the central pair. To keep this paper close to being self-contained, 
we provide  most proofs. 

The polynomials spaces $\calP(X)$ and $\calD(X)$ are best described
in terms of a partition of the power set $2^X$ 
into two disjoint sets of long subsets $L(X)\subset 2^X$ and short subsets
$S(X)=2^X\bks L(X)$: 
\begin{eqnarray*} L(X) & \eqbd &  \Set{ Y\subset X \,:\,
  Y\cap B \neq \emptyset, \quad \forall B \in\B(X) }, \\
 S(X) & \eqbd & \Set{Y\subset X \,:\,
\rank(X\backs Y) =n }. \end{eqnarray*}
Note that the elements of $S(X)$ are exactly the independent 
sets of the matroid dual to $X$, and those of $L(X)$ are its
dependent sets, as the independence of a set $Y$ in the dual
matroid is equivalent to the set $X\setminus Y$ being of full
rank in the original matroid. 
Using the notation~(\ref{def:p_Y}) for polynomials $p_Y$, 
the spaces $\calP(X)$ and $\calD(X)$ are defined as follows:
\begin{eqnarray*} \calD (X) & \eqbd & \Set{ f\in \Pi \,:\,
p_Y(D)f = 0, \quad \forall Y\in L(X)},   \\
\calP (X) & \eqbd &  {\rm span} \Set{ p_Y \,:\, Y \in S(X)}.
\end{eqnarray*}
Immediately, Theorem \ref{thm_db'} and Corollary~\ref{cor_db'} imply as a 
special case the following two statements:

\begin{theorem} \label{thm_pi_d}
Let $V(X,\lam)$ denote the vertices of a generic 
$X$-hyperplane arrangement associated with $X$. Then 
$\Pi(V(X,\lam)) \subseteq \calD(X)$.
\eop \end{theorem}

\begin{corollary} \label{cor_lb_d} $\dim \calD(X) \geq \# \B(X)$. 
\eop \end{corollary}

We now let $\calJ(X)$ denote the ideal generated by the long
polynomials: $$\calJ(X)\eqbd \Ideal \Set{ p_Y \,:\, Y \in L(X)}.$$
Then $\calD(X)$ could also be defined as the polynomial kernel of $\calJ(X)$:
$\calD(X)=\Pi\cap \ker\calJ(X).$
It is easy to check that the only common zero of the
long polynomials is $0$, i.e., that $\var(\calJ(X)) = \Set{0}$. 
Since the ideal $\calJ(X)$ is homogeneous, this  tells us  
(see Section~\ref{sec_ideals}) that the kernel of $\calJ(X)$ is 
finite-dimensional and is a subspace 
of $\Pi$. Thus,
$\calD(X)$ is precisely the kernel of $\calJ(X)$:
 \begin{equation}  \ker \calJ(X) = \calD (X) \subset \Pi. \label{D=ker(J)}
\end{equation}

\begin{theorem} \label{thm_sum} $\calP(X)+\calJ(X)=\Pi$. 
\eop \end{theorem}

\nt{\bf Proof, \cite{DR_TAMS}.\/} Set $F\eqbd \calP(X)+\calJ(X)$, and assume
that $F$ is proper in $\Pi$. Since $F$ is homogeneous,  its orthogonal complement $F{\perp}$ in $\Pi$ (with respect to the pairing
(\ref{inner_prod})) contains non-zero polynomials. We will show that
$F{\perp}$ is $D$-invariant, i.e., closed under differentiations.
This will imply that $F{\perp}$ must contain the constants,
which is absurd, since $\calP(X)$, hence $F$, contains the constants.

We need thus to prove that, for $p\in F{\perp}$ and $a\in \Rn$,
$D_ap\perp F$. First, since $$p\in F{\perp}\, \subset \,\calJ(X){\perp}=\ker \calJ(X),$$
it follows that $D_ap\perp\calJ(X)$ (since $\calJ(X)$ is an ideal
and the kernel of an ideal is always $D$-invariant).
It remains to show that $D_ap\perp\calP(X)$, i.e., that
$D_ap\perp p_Y$, for every short $Y$. Fix such $Y$ and choose
$B\in \B(X\bks Y)$. Since we can write
$a=\sum_{b\in B}c(b)b$, it suffices to prove that
$D_bp\perp p_Y$ (for every $b\in B$). 
Now,
$$\inpro{D_bp,p_Y}=\inpro{p,p_bp_Y}=\inpro{p,p_{Y\cup b}}.$$
Since $Y\cup b\subset X$, and since $F$ contains, by assumption,
every polynomial $p_{Y'}$, $Y'\subset X$, it follows that
$p_{Y\cup b}\in F$, hence $\inpro{p,p_{Y\cup b}}=0$.
\eop


\begin{corollary} \label{cor_lb_p} $\dim \calP(X) \ge \dim \calD(X)$.
\end{corollary}

\nt {\bf Proof.\/} This inequality follows by applying Result~\ref{res_sum} 
to the sum from Theorem~\ref{thm_sum} and by recalling the
formula~(\ref{D=ker(J)}) that identifies $\ker\calJ(X)$ with $\calD(X)$.
  \eop

Next, consider the set of {\it facet hyperplanes\/} of $X$:
\begin{equation}\label{defdualhyp}
\calF(X)  \eqbd \Set{ H  \,:\, H \mbox{ is a subspace of } \R^n, \ 
\dim H = n-1, \ \span(X\cap H) =H}. 
\end{equation}
Note that these hyperplanes are in fact parallel to the facets of the zonotope 
$Z(X)$, which explains this terminology.  Given any facet hyperplane  
$H \in \calF(X)$, let $\eta_H$ be a non-zero normal to $H$: 
$\eta_H \perp H$. We also define
\begin{equation}\label{defmh}
m(H) \eqbd  m_X(H) \eqbd  \#(X\backs H).
\end{equation}
Define 
$$\calI(X) \eqbd  \Ideal \Set{p^{m(H)}_{\eta_H} \,:\, H\in \calF(X)},$$
where, as above, $p_x: t\mapsto x\cdot t$. 
Then we have the following theorem:

\begin{theorem} \label{thm_p_i}
$\calP(X) \subseteq \ker \calI(X)$. 
\end{theorem}

\nt {\bf Proof.\/}
We only need to check that
$$ D^{m(H)}_{\eta_H} (p_Y)=0 \quad \forall H \in \calF(X) \quad 
\forall Y\in S(X), $$
i.e., that any generator of $\calI(X)$ acting as a differential operator
annihilates any generator of $\calP(X)$. Indeed, note that, for
a single vector $\xi\in X$, we get $D_{\eta_H}(p_\xi)=0$ if and only if 
$\xi\perp
\eta_H$ or, equivalently, $\xi\in H$. So, a polynomial of the form $p_Y$ can survive
$m(H)$ differentiations in the direction $\eta_H$ only if $\# (Y\backs H)\geq
\# (X \backs H)$, i.e., only if $\# (Y\backs H)= \# (X \backs H)$ 
since $Y\subseteq X$. The last equality is possible if and only if
all vectors in $X\backs Y$ belong to $H$.  But $Y$ is a short subset of $X$,
so its complement is of full rank. Contradiction! Hence $p_Y$ is annihilated
by $D_{\eta_H}^{m(H)}$.   \eop

\begin{corollary} \label{cor_ub_p} $\dim \calP(X) \leq \dim \ker \calI(X)$. 
\eop \end{corollary}

\begin{theorem}[\cite{BDR_PJM}] \label{thm_dim_i} $\dim \ker \calI(X) \leq 
\# \B(X)$. 
\end{theorem}

\nt{\bf Sketch of proof.} The proof is by induction on $\#X$ and $n$. 
Assuming that this statement
is correct for $X$, we define $X'\eqbd X\cup\{\xi\}$, and consider
for every facet hyperplane  $H\in \calF(X)$
the space $P_H\eqbd \calP((X\cap H)\cup\{\xi\})$. If $\xi\in H$, then
$P_H=0$; otherwise, $P_H$ has positive dimension. Note that each
$B\in\B(X')$ lies either in $\B(X)$ (in case it does not contain 
$\xi$), or else in a unique $(X\cap H)\cup\{\xi\}$, $H\in \calF(X)$.
Therefore,
$$\dim \calP(X)+\sum_{H\in \calF(X)}\dim P_H=
\#\B(X').$$
We then define a map $T$ as follows:
$$T:\ker I(X')\to\times_{H\in\calF(X)}P_H  \;\; : \;\;
f\mapsto (D^{m(H)}_{\eta_H}f)_{H\in \calF(X)}.$$
The kernel of this map is, directly from the definition, $\ker I(X)$, hence, 
by induction, $\calP(X)$. Our previous computation then shows that
$$\dim \ker I(X')\le \dim\ker T+\dim\ran T=\#\B(X').$$
The only missing item in the argument is to show that the map $T$ is
well-defined, i.e., that $D_{\eta_H}^{m(H)}\ker I(X')\subset
P_H$. This is trivially true in case
$\xi\in H$. Proving the above for the case $\xi\in H$ is the hard part
of the proof, which is omitted here. Some of these missing details
are discussed as a part of the proof of Theorem~\ref{thm_dim_i+} in the next section. See \cite{BDR_PJM} for  details.   \eop

We are now in a position to prove the main theorem of this section.

\begin{theorem} \label{thm_main} \hfill
\bd
\item{(1)}\quad 
$\dim \calP(X) = \dim \calD(X) = \#\B(X)$.
\item{(2)}\quad
The map $p \mapsto \inpro{p, \cdot}$ is a bijection between
$\calD(X)$ and $\calP(X)'$. 
\item{(3)}\quad
$\calD(X) = \Pi(V(X,\lam))=\ker \calJ(X)$. 
\item{(4)}\quad
The point set $V(X,\lam)$ is correct for $\calD(X)$ as well as for
$\calP(X)$.
\item{(5)}\quad
$\calP(X) = \ker \calI(X)$.
\item{(6)}\quad 
$\calP(X)\oplus\calJ(X)=\Pi$.
\ed
\end{theorem}

\nt{\bf Proof.\/} Putting together the inequalities obtained 
in Corollaries~\ref{cor_lb_d}, \ref{cor_lb_p}, \ref{cor_ub_p}
and in Theorem~\ref{thm_dim_i}, we get
$$ \# \B(X)\leq \dim \calD(X) \leq \dim \calP(X) \leq \dim \ker\calI(X) 
\leq \# \B(X).   $$ 
This shows that equalities must hold throughout.
Invoking Theorems~\ref{thm_pi_d},~\ref{thm_sum} and~\ref{thm_p_i}, along with 
 Result~\ref{res_pi(sig)}, Corollary~\ref{cor_duala} and Result~\ref{res_sum}, we obtain the remaining  claims of this theorem. \eop

\nt {\bf Remark.\/} Let us assume that $X=B$, with $B$ some basis for
$\Rn$. Then $S(X)=\{\emptyset\}$, hence $\calP(X)=\span\{1\}=\Pi_0^0$.
On the other hand, $\{b\}\in L(X)$, for every $b\in B$, and hence
$\calJ(X)$ contains the linear polynomials $\{p_b: b\in B\}$. Thus,
$\calJ(X)$ is the maximal ideal $\{p\in \Pi: p(0)=0\}$, and the
decomposition $\calP(X)\bigoplus \calJ(X)=\Pi$ becomes obvious.
There are other cases when this decomposition can be obtained directly:
for example, when $X$ is in general position. However, for a general
$X$, this decomposition is non-trivial.

\medskip
As a by-product of Theorem~\ref{thm_main}, we obtain an additional result
that characterizes the least space obtained from integer points in
the half-open half-closed zonotope in the case when $X$ is unimodular,
i.e., when all vectors in $X$ have only integer components and every
basis of $X$ is invertible over $\Z$ (see Section~\ref{sec_zono}).  
We recall that the zonotope $Z(X)$ is defined as the image 
$Z(X)\eqbd  X([0,1]^X)$ of the unit cube $[0,1]^X$ under the map 
$X: \R^X\to \R^n$. 

Thus assume that $X$ is unimodular and consider its zonotope and 
associated hyperplane arrangements. In the context of  hyperplane 
arrangements, a set of interest is the
vertex set $V(X,\lam)$ of the arrangement, whose precise geometry
depends on the vector $\lam\in \C^X$. For a generic $\lam$, the vertex set
$V(X,\lam)$ is of {\it maximal\/} cardinality $\#\B(X)$.
The dual vertex set, $\calZ(X,t)$, is parameterized by $t\in \R^n$.
As we will shortly see, for a generic $t$ this set is of {\it minimal\/} 
cardinality $\#\B(X)$. Let us begin with a definition:
$$\calZ(X,t)\eqbd \{\alp\in \Z^n:\ t-\alp\in Z(X)\}.$$
We consider $t$ to be {\it generic\/} if it does not lie in any of the 
hyperplanes $$\Z^n+H,\quad H\in \calF(X).$$
If $t$ is generic, then it is well known (see, e.g.,~\cite{BHR})
that $$\# \calZ(X,t) = {\rm Vol} (Z(X))=\#\B(X).$$
We will assume that $t$ is fixed and generic and will occasionally 
denote  $\calZ(X,t)$ simply by $\calZ(X)$.

\begin{theorem} [\cite{BDR_PJM}] \label{thm_intpoints} Let $X$ be unimodular. Then 
$\Pi(\calZ(X))=\calP(X)=\ker \calI(X)$, hence $\calZ(X)$ is correct for
$\calP(X)$ as well as for $\calD(X)$.
\end{theorem} 

\nt{\bf Proof.\/} With $\calZ(X)=\calZ(X,t)$, we already know that 
$ \# \calZ(X) = {\rm Vol} (Z(X)) = \#\B(X). $
Recall that for a given $\sig \subset \R^n$,  we defined
$$\Exp(\sig) \least \bdeq \Pi (\sig).$$
By Theorem~\ref{res_pi(sig)}, $\dim\Pi(\sig)=\#\sig$ for any set $\sig$, 
so, for a unimodular $X$, we get
$$\dim\Pi(\calZ(X))=\#\calZ(X)=\#\B(X).$$
Since both spaces $\Pi(\calZ(X))$ and $\calP(X)$
have the same dimension, we only need to 
prove that one is included in the other.
We will show that $\Pi(\calZ(X)) \subset \ker \calI(X)$.
To this end, we recall Corollary~\ref{cor_annihil}:
if  $f,g\in \Pi$ satisfy  
$f_\uparrow = g$, and if $f|_\sig=0$, then 
$$g(D)(\Pi(\sig))=0.$$
We choose $g$ to be one of the generators of $\calI(X)$, i.e.,
$$g: t\mapsto (\eta_H\cdot t)^{m(H)}=p_{\eta_H}^{m(H)}(t),\quad H\in \calF(X).$$
We need to find $f$ such that $f\most=g$ and $f$ vanishes on $\calZ(X)$.
Once we manage to do so for every $g$ as above, we are done.
To this end, we fix $H \in \calF(X)$, and will define $f$ as 
(with $\eta\eqbd \eta_H$) 
$f \eqbd  (p_\eta +c_1)(p_\eta +c_2)\cdots (p_\eta +c_{m(H)})$,
with $(c_i)_i$ some constants. Obviously, for such $f$ we always
have that $f\most =g$. We need also to ensure that $f$ vanishes on
$\calZ(X)$. In the argument below we assume for convenience that
all vectors $X\bks H\bdeq \{y_1,\ldots,y_{m(H)}\}$ lie on one side of $H$. 
It is then straightforward to see that the zonotope $Z(X)$ lies between
the hyperplane $H$, and the hyperplane
$$H'\eqbd H+\sum_{i=1}^{m(H)}y_i.$$
A simple consequence of the unimodularity is that there are exactly
$m(H)-1$ translates of $H$ that lie properly between $H$ and $H'$ and
contain integers. Precisely, these are the hyperplanes
$$H+\sum_{i=1}^jy_i\bdeq H+c_j,\quad j=1,\ldots,m(H)-1.$$
Since $t$ is generic, it does not lie on any of these hyperplanes, hence we
may assume without loss that it lies between $H$ and $H+y_1$. We
conclude that
$$\calZ(X,t)\subset \cup_{j=0}^{m(H)-1} (H+c_j),\quad c_0\eqbd 0,$$
hence that the polynomial
$$f\eqbd \prod_{j=0}^{m(H)-1}(p_\eta+c_j)$$
vanishes on $\calZ(X,t)=\calZ(X)$. This proves that
$\Pi(\calZ(X))=\calP(X)=\ker \calI(X)$, the second equality by
Part~(5) of Theorem~\ref{thm_main}.

The correctness of $\calZ(X,t)$ for $\calP(X)$ follows then from
Result~\ref{res_pi(sig)}, while its correctness for $\calD(X)$ 
follows from the duality between $\calP(X)$ and $\calD(X)$ 
proved in Part~(2) of Theorem~\ref{thm_main} and from Result~\ref{cor_duala}.
\eop

\nt{\bf Additional historical remarks.\/}
The space $\calD(X)$ was introduced in~\cite{BH_JdAM}.
The inequality  $\dim\calD(X)\ge \#\B(X)$ was first proved by Dahmen 
and Micchelli in~\cite{DM_Studia} by induction  on $\#X$ and on $n$. 
A non-inductive analytic argument is given in~\cite{BenR_TAMS}.  
The equality~$\dim\calD(X)=\#\B(X)$ is also due to Dahmen and 
Micchelli~\cite{DM_Studia}.  They subsequently provided, 
in~\cite{DM_Adv}, a very elegant proof for the inequality 
$\dim\calD(X)\le\#\B(X)$, which uses the matroidal structure of $X$. 
In~\cite{DR_TAMS}, the
space $\calP(X)$ is proved to be dual to every space of the form
$\calD(X,\lam)$, with the definition of the latter obtained from the
definition of $\calD(X)$ by replacing each $p_x$ by $p_{x,\lam}$; here
$\lam$ need not be generic. The space $\calD(X,\lam)$ plays an important
role in the theory of {\it exponential box splines,\/} \cite{R_CA_1988},
but not in this paper.

\subsection{\label{sec_basis}Homogeneous basis and Hilbert series
for $\calP(X)$}

Let $\Pi^0_j$ be the space of homogeneous polynomials of degree $j$
(in $n$ variables).
Since both $\calP(X)$ and $\calD(X)$ are {\it graded\/} or 
{\it homogeneous,\/}
i.e., are spanned by homogeneous polynomials and since the 
pairing~(\ref{inner_prod}) respects grading, the isomorphism 
shown in part (2) of Theorem~\ref{thm_main} implies that, for every $j$,
$$\dim(\Pi_j^0 \cap \calD(X)) =\dim(\Pi_j^0 \cap \calP(X)).$$
We refer to the homogeneous dimensions of the space $\calP(X)$ as
the {\it central Hilbert  series of $X$:\/} 
$$h_X: \N \to \Z_+ \,:\, j \mapsto \dim(\Pi_j^0 \cap \calP(X)).  $$
Note that  
$$\sum_j h_X(j) = \# \B(X).$$
The adjective ``central'' is chosen in anticipation of the introduction 
of two other Hilbert series that will be labeled ``internal'' and 
``external'' respectively.

We focus now on building a homogeneous basis for $\calP(X)$, which will enable
us to compute the homogeneous dimensions $h_X(j)$ of $\calP(X)$. 
We will see soon that $h_X$ can be computed directly by studying
the dependence/independence relations among the vectors in $X$.

\begin{example} Let
$$X=\left[ \begin{array}{rrrrrr}1&0&0&1&1&0\\ 0&1&0&-1&0&1\\ 0&0&1&0&-1&-1
\end{array}
\right].$$
Then $\#\B(X)=16$, and 
$$h_X=(1,3,6,6).$$
\qed
\end{example}

\nt  {\bf An algorithm for computing $h_X$.\/}
First, we impose an arbitrary order $\prec$ on $X$. Then we associate each 
$B \in \B(X)$ with the  homogeneous polynomial $p_{X(B)}$, where 
\begin{equation}\label{defxb}
X(B)\eqbd  \Set{ y\in X \,:\, y \notin \span\Set{b\in B \,:\, b \preceq y}}.
\end{equation}
Note that $X(B) \in S(X)$, since $B\subset (X\bks X(B))$.
Define
$$\val B\eqbd \#X(B).$$
Then
$$h_X(j)=\#\{B\in\B(X): \val B=j\}.$$
This assertion follows from the stronger assertion in Theorem~\ref{thm_basis}
below. The algorithm, incidentally, draws an intimate relation between the 
Tutte  polynomial of $X$~\cite{Biggs} and the central Hilbert function.

\begin{example} \quad Let
$X\eqbd \left[ \matrix{1 &\ 1&\ 0& \ 1  \cr  0&\ 0&\ 1&\ 1 \cr} \right] \bdeq 
 [x_1\  x_2\  x_3\  x_4]$ be an ordered multiset. Then 
 $$ X([x_1,x_3]) = \emptyset, \qquad 
    X([x_2,x_4]) = \Set{x_1, \ x_3}.$$
 The algorithm easily produces the Hilbert series
$$h_X=(1,2,2,0,0,\ldots)$$
 \qed
\end{example}

\begin{theorem} [\cite{DR_TAMS}] \label{thm_basis} \quad
The set  
\begin{equation}
\Set{ Q_B\eqbd p_{X(B)}\,:\, B \in \B(X)} \label{pYB}
\end{equation}
is a  basis for $\calP(X)$.
\end{theorem}

\nt{\bf Proof.\/} It is clear that  
$p_{X(B)}\in \calP(X),$
since $X(B)\cap B = \emptyset$, i.e., $X(B) \in S(X)$.
Since $\dim\calP(X)=\#\B(X)$, it is sufficient to 
show that 
$\Set{ Q_B \,:\, B \in \B(X)}$ is linearly independent. 
We will prove this by induction on $\# X$, with the case $\#X=n$ being
trivial.

Assume, thus, the set in (\ref{pYB}) to be linearly independent.
Let $X' = X\cup \Set{\xi}$ where $\xi$ is the last element in $X'$.
The induction step requires us to show that, given $p\in \calP(X)$
and $B_0\in {\B(X')\backs \B(X)}$,
if, for some constants $(a(B))_B$,
$$p + \sum_{B\in \B(X')\backs \B(X)} a(B)\, p_{X(B)} =0,$$
then $a(B_0)=0$. To this end,
we define
 $$B_0 \bdeq B_0' \cup \Set{\xi}, \quad 
 H \eqbd  \span\, B_0'\in\calF(X).$$
Let $\eta$ be a non-zero vector such that $\eta \perp H$. Then
\begin{equation}
0=D^{m(H)}_\eta 0 = 
D^{m(H)}_\eta p + \sum_{B\in \B(X')\backs \B(X)} a(B) D^{m(H)}_\eta p_{X(B)},
\label{Deq}
\end{equation}
where $m(H) \eqbd  \#(X\backs H)$.  Since  $D_\eta^{m(H)}\calP(X)=0$ 
by~Theorem \ref{thm_p_i}, we conclude that  $ D^{m(H)}_\eta p =0.$
Next, consider 
$$D^{m(H)}_\eta p_{X(B)}, \quad B\in \B(X')\backs \B(X).$$
We know  that $\#(X' \backs H) = m(H) +1$. 
If $\#(B\cap H) < n-1$, then $X'\backs X(B)$ contains at least 2 vectors from 
$X' \backs H$. 
So, $ \#( X(B) \cap (X'\backs H)) < m(H)$. Consequently,
$$ D^{m(H)}_\eta p_{X(B)} =0.$$
If $\#(B\cap H) = n-1$, then $ \#(X(B) \backs H) = m(H)$ so that 
$$  D^{m(H)}_\eta p_{X(B)} = c(B) \, p_{X(B)\cap H},$$
for some non-zero coefficient $c(B)$. From equation (\ref{Deq}), we get 
$$0= \sum_{B\in K } a(B) c(B) \, p_{X(B)\cap H},$$
where 
$K\eqbd  \Set{B \in \B(X')\backs \B(X) \,:\, \#(X(B) \backs H) =m(H)}$. 
Now, it is easily observed that, with $W\eqbd \xi\cup(X'\cap H)$, the 
polynomials
$$p_{X(B)\cap H},\quad B\in K,$$
are exactly the polynomials in the homogeneous basis for $\calP(W)$,
with the order on $W$ being the order induced from $X$. Our induction
hypothesis implies (since $\#W<\#X'$) that those polynomials are independent,
hence that $a(B_0)=0$.  \eop

\nt {\bf Remark.\/} The construction provides us with the direct sum
decomposition
$$\calP(X\cup\{\xi\})=\calP(X)\bigoplus (\bigoplus_{\xi\not\in H\in
\calF(X)}\ p_{X\bks H}\calP((X\cap H)\cup\{\xi\}).$$
The decomposition allows us to compute
the Hilbert series $h_{X'}$ by summing the Hilbert series of the various
summands:
$$h_{X\cup\xi}(j)=h_X(j)+\sum_{\xi\not\in
H\in\calF(X)}m_X(H)+h_{X_H}(j),$$
with $X_H\eqbd (X\cap H)\cup\xi$.
This means that we do not need to impose one fixed order
on $X$: once $\xi$ is known to be placed last, the order of the remaining
elements can be chosen separately (hence independently) for each summand. 
This is consistent with the known invariance of the Tutte polynomial,
\cite{Biggs}.

\section{\label{sec_ext}External algebra}

\subsection{\label{sec_main+}Main results}

Recall that, in Section~\ref{sec_linalg}, we let $\I(X)$ denote 
the collection of all independent subsets of $X$ and let 
 $B_0\subset \R^n$ be a fixed ordered basis. We also denoted
$$X' \eqbd X \cup B_0$$ and defined a bijection 
$\ex \,:\, \I(X) \to \B_+(X)\subset \B(X')$
to the set $\B_+(X)$ of all external bases of $X$. We now define
\begin{eqnarray*} L_+(X) & \eqbd & \Set{Y\subset X' \,:\, 
Y \cap B \neq  \emptyset, \quad B \in \B_+(X)}, \\
\calP_+(X) & \eqbd & \span \Set{p_Y \,:\, Y\subset X}, \\
\calD_+(X) & \eqbd & \Set{f \in \Pi \,:\, p_Y(D)f=0, \
\forall Y \in L_+(X)}.
\end{eqnarray*}
Our goal is to show that $\calD_+(X)$ and $\calP_+(X)$ are dual to 
each other, and to determine their annihilating ideals. The main
result will be established in several steps, analogously to 
Section~\ref{main_central}.

Let $V(X',\lam)$ be the set of vertices of a generic $X'$-hyperplane 
arrangement.  With a slight abuse of notation we denote by $V_+(X,\lam)$ the 
subset of $V(X',\lam)$ that corresponds to $\B_+(X)\subset \B(X')$.
Applying Theorem~\ref{thm_db'} and its Corollary~\ref{cor_db'}
to this case (i.e., with $X$ there replaced by $X'$ here, and $\B'$
there being our $\B_+(X)$) we obtain the following results.

\begin{theorem}  \label{thm_pi_d+} $\Pi(V_+(X,\lam)) \subseteq \calD_+(X)$.   \eop
\end{theorem}

\begin{corollary} \label{cor_lb_d+}
$\dim \calD_+(X) \geq \# \B_+(X)$. \eop
\end{corollary}

Define 
$$\calJ_+(X) \eqbd \Ideal \Set{p_Y \,:\, Y\in L_+(X)}.$$
Note that, (almost) directly from the definition of $\calD_+(X)$,
$\ker \calJ_+(X) = \calD_+(X)$.

\begin{theorem} \label{thm_sum+} $\calP_+(X)+\calJ_+(X)=\Pi$. 
\end{theorem}

\nt {\bf Proof.\/} We start with the fact that
$\calP(X) +  \calJ(X) = \Pi$, which is established in Theorem~\ref{thm_sum}.
Since $\calP_+(X) \supset \calP(X)$, we conclude that
$$ \calP_+(X) +  \calJ(X) = \Pi.$$
So, we need to prove  that
$$\calJ(X) \subset \calP_+(X) +  \calJ_+(X).$$
Let $Y\in L(X), \ f\in \Pi$. Since every polynomial
in $\calJ(X)$ is a combination of polynomials of the form
$p_Yf$, it  suffices to prove
that 
$$p_Yf \in  \calP_+(X) +  \calJ_+(X),$$ 
a claim that we prove by
reverse induction on $\#Y$. Thus, assume that the claim is correct for every 
$Y' \in L(X)$ such that $\#Y' > \#Y$. 
Put $S \eqbd  \span(X \backs Y)$. Then $\dim S < n$, since $Y$ is long. 
Let $I \subset X \backs Y$ be a basis for $S$ and 
$B \eqbd  \ex(I)$. Since $B$ is a basis for $\R^n$, 
$$\Ideal \Set{p_b \,:\, b \in B} + \Pi_0^0 = \Pi.$$
So, we can write $f$ in the following form:
$$f = c_0 + \sum_{b\in B} p_b f_b, \quad f_b \in \Pi.$$
Consequently,
$$ p_Yf = c_0 p_Y + \sum_{b\in B} p_{Y\cup \Set{b}} f_b. $$
We claim that each term above belongs to  
$\calP_+(X) +  \calJ_+(X) $. 
Since $Y\subset X$, it is clear that $p_Y \in \calP_+(X)$.
Now, for $p_{Y\cup \Set{b}}$, we have 
either $b \in I$ or $b\in B_0$.

{\bf Case I.\/} If $b \in I \subset X$, then 
$Y'\eqbd  Y \cup \Set{b} \subset X$.
By induction,
$$p_{Y\cup \Set{b}} f_b \in \calP_+(X) +  \calJ_+(X). $$

{\bf Case II.\/} Let $b\in B_0$. We show that
$$p_{Y\cup \Set{b}} f_b \in \calJ_+(X), $$
and to this end it is enough to show that 
$Y \cup \Set{b} \in L_+(X)$. 
Let $B'\in \B_+(X)$. If $Y \cap B' = \emptyset$, then 
$B'\cap X\subset X\bks Y \subset S $. Hence $B'=\ex(I')$,
for $I'\subset S$. Thus $\span I'\subset \span I$, and the definition
of the extension map implies that in such case we always have that
$B_0\cap \ex(I)\subset B_0\cap \ex(I')$. Consequently, $b \in B'$,
and thus
$$ \gp{ Y\cup \Set{b} } \cap B' \neq \emptyset.$$ 
We conclude that $Y\cup\Set{b}\in L_+(X)$, hence that, directly from
the definition of $\calJ_+(X)$,
$p_{Y\cup\Set{b}}f_b\in \calJ_+(X)$; {\it a fortiori\/} 
$p_{Y\cup\Set{b}}f_b\in \calP_+(X)+\calJ_+(X)$.
 \eop
 
\nt{\bf Remark.\/} Note that the only property of the extension $\ex$
that was used is that once $\span I'\subset \span I$, then $B_0\cap I\subset 
B_0\cap I'$. It is probably easy to show that every extension of such 
type is a greedy extension with respect to some ordering of $B_0$.
\qed
 
Invoking Result \ref{res_sum}, we obtain

\begin{corollary} \label{cor_lb_p+}  $\dim \calP_+(X) \ge \dim \calD_+(X)$.
\eop \end{corollary}

The corollary implies that $\dim\calP_+(X)\ge \#\B_+(X)$.
This last estimate can be proved directly: Order $X'$ so that $B_0$ 
is placed after $X$, and the internal order within $B_0$ is retained. 
Then follow the construction of a homogeneous basis for $\calP(X')$
from Section~\ref{sec_basis}.
Observe that a polynomial in that basis
is of the form $p_{X'(B)}$, $X'(B)\subset X'$, and that $X'(B)$ is then a 
subset of $X$ if and only if $B\in \B_+(X)$. Thus
$$\Set{p_{X'(B)}:\ B\in \B_+(X)}\subset \calP_+(X),$$
and we get the desired bound from the linear independence of these
polynomials. We will come back to this issue later,
since the polynomials above form a {\it basis\/} for $\calP_+(X)$,
and we will use the cardinality of the sets
$X'(B)$, $B\in \B_+(X)$, in order to provide an algorithm for computing
the forthcoming external Hilbert series $h_{X,+}$ of $X$.

\bigskip
Recalling (\ref{defdualhyp}) and (\ref{defmh}), we define
$$\calI_+(X) \eqbd \Ideal \Set{p_{\eta_H}^{m(H)+1} \,:\, 
H \in \calF(X),\ \eta_H\perp H }.$$

\begin{theorem} \label{thm_p_i+} $\calP_+(X) \subseteq \ker \calI_+(X)$. 
\end{theorem}

\nt{\bf Proof.\/} Given any $Y\subset X$ and any facet hyperplane $H$,
we have that $D_{\eta_H}(p_Y)=p_{Y\cap H} D_{\eta_H} p_{Y\bks H}$.
The result then follows from the fact that $\#(Y\bks H)\le m(H)$.
\eop

\begin{corollary} \label{cor_ub_p+}
$\dim P_+(X) \leq \dim \ker \calI_+(X)$. 
\eop \end{corollary}

\begin{theorem} \label{thm_dim_i+} $\dim \ker \calI_+(X) \leq \# \B_+(X)$. 
\end{theorem}

\nt 
Before we embark on a proof of this theorem, we must make a few auxiliary
statements first. For our next result, we will use the symbol
$$ \calP_N^0(S)\eqbd \spa\{p_Y : Y\subset S , \#Y=N  \}$$
to denote the space of homogeneous polynomials of degree $N$ in the variables 
$S$,  for any (possibly infinite) set $S$ and any nonnegative integer $N$.

\begin{proposition} \label{prop_hom_total} 
Let $I$ be an ideal of $\Pi$ and let $S$
be a subspace of $\R^n$ of dimension $d\geq 2$. 
Let $S_1$, $\ldots$, $S_k$ be distinct subspaces of $S$, 
each of dimension $d-1$. Suppose that, for $n_1$, $\ldots$, 
$n_k \in \N$, the ideal $I$ contains all homogeneous polynomials 
in variables  $S_i$ of degree $n_i$:
$$  \calP_{n_i}^0 (S_i) \subset I.  $$
Then $$ \calP_N^0 (S) \subset I 
\quad  {\rm whenever} \quad 
(N+1)(k-1) \geq   \sum_{i=1}^k n_i. $$
\end{proposition}

\nt{\bf Proof.\/} Pick $l\in S$. We need to show that the polynomial 
$p_l^N$ lies in the ideal $I$. Choose a subspace $V\subset S$ of 
dimension $2$  such that $l\in V $ but $V\not\subset S_i$ for all 
$i=1, \ldots, k$.
The dimension formula $$\dim V+\dim S_i = \dim (V\cap S_i) + \dim (V+S_i)$$
immediately implies that $\dim (V\cap S_i)=1$ for each $i=1, \ldots, k$.
So, for each $i$, there exists a nonzero vector $h_i \in V\cap S_i$.
By the assumption of the Proposition, 
$$ p_{h_i}^{n_i} \in I, \qquad i=1, \ldots, k.$$
Observe that
$$ H_i\eqbd p_{h_i}^{n_i} \calP_{N-n_i}^0(V) \subset I \cap \calP_N^0(V), 
\qquad i=1, \ldots, k.$$   
We will now argue that $\sum_{i=1}^k H_i = \calP_N^0 (V)$. 
Indeed, if not, then there exists a nonzero polynomial 
$q\in \calP_N^0(V)$ such that $\langle q, p \rangle =0$ 
for all $p\in H_i$. This means 
\begin{equation} 0 = \langle q, p_{h_i}^{n_i} p \rangle=
\langle p(D) q , p_{h_i}^{n_i} \rangle 
\label{pq_eval} \end{equation}
for any $p\in \calP_{N-n_i}^0(V)$.   
But the last expression in~(\ref{pq_eval}) is, up to the factor $n_i!$,
equal to $(p(D)q)(h_i)$. Our last setup can therefore be reformulated as 
a univariate problem: there exists a nonzero polynomial 
$q\in \calP_N(\R)$ and distinct points $h_i\in \R$, $i=1, \ldots, k$ such that
\begin{equation} 
 (p(D)q)(h_i)=0, \qquad i=1, \ldots, k  \quad \hbox{\rm for all} 
\quad p\in \calP_{N-n_i}(\R).
\label{pq_system}
\end{equation}
But this is a Vandermonde linear homogeneous system of 
$\sum_{i=1}^k (N-n_i+1)$ equations  in $N+1$ unknown coefficients 
of $q$, which has a nontrivial solution if and only if 
$$ \sum_{i=1}^k (N-n_i+1) < N+1 \quad {\rm iff} \quad 
  Nk-\sum_{i=1}^k n_i +k < N+1 \quad {\rm iff} \quad
(N+1)(k-1) < \sum_{i=1}^k n_i, $$
contrary to the assumption of this Proposition. Thus 
$p_l^N \in I$ for all nonzero vectors $l\in S$, hence every homogeneous
polynomial in $S$ of degree $N$ is in $I$ and therefore 
$\calP_N^0(S)\subset I$.
    \eop

\begin{corollary}  \label{cor_hom_i+}
 Let $Y\subset X$, and let $0\neq \eta \perp \span Y$. 
Then $p_\eta^{\# (X\backs \span Y)+1 } \in \calI_+(X)$.   
\end{corollary}

\nt{\bf Proof.\/} We run the proof by induction on $n-\dim (\spam Y)$.
When $Y$ spans  a hyperplane, we have $p_\eta^{\# (X\backs \spam Y)+1} 
\in \calI_+(X)$ by definition of $\calI_+(X)$.
When $\dim (\spam Y) \leq n-2$, we denote $S\eqbd (\spam Y){\perp}$ 
and consider all possible ways to  add one more vector to the set $Y$ 
to increase the dimension of $\spam Y$. Call the orthogonal complements
of the spans of these sets $S_1$ through $S_k$. Note that addition
of different vectors $x$ to $Y$ may produce the same subspace 
$\span (Y\cup \{x\})$  and therefore the same orthogonal complement. 
If that  is the case, we list such a subspace $S_i$ only once.

With $\eta \in S$ and $S_1$ through $S_k$ subspaces of $S$ of
codimension $1$, we are now in the setting of  
Proposition~\ref{prop_hom_total}, so may conclude that
$$ p_\eta^N \in \calI_+(X) $$
whenever $(N+1)(k-1)\geq \sum_{i=1}^k  n_i$, where 
$n_i = \# (X\backs (S_i{\perp}))+1$ by the inductive hypothesis. 
Note that the count $ \#(X\backs (S_i{\perp}))$ performed for all $i$ 
accounts for every vector of $(X\backs (S{\perp}))=X\backs\span Y$ exactly 
$k-1$ times, hence 
\begin{eqnarray*}
 \sum_{i=1}^k  \#(X\backs (S_i{\perp})) & = &  (k-1)\cdot 
 \# (X\backs \span Y), \quad {\rm hence} \\
  (k-1)(N+1) & \geq & (k-1)\cdot  \# (X\backs \span Y) + k , 
 \quad \hbox{\rm or, equivalently,} \\
N &\geq & \# (X \backs \span Y) +1/(k-1). 
  \end{eqnarray*}
The last inequality is satisfied whenever $N\geq \# (X\backs \span Y)+1$,
so we are done.   \eop

We are now in a position to give a proof of Theorem~\ref{thm_dim_i+}.
\medskip

\nt{\bf Proof of Theorem~\ref{thm_dim_i+}.\/}  {\bf Step I.\/} We append to $X$ an auxiliary basis $B_0$, and obtain
$X'\eqbd X\cup B_0$. We choose $B_0$ to be in general position with
respect to $X$.
Note that  $$\calI(X') \subset \calI_+(X).$$
Indeed, let $H\in \calF(X')$ be a facet hyperplane of $X'$
and let $\eta$ be normal to $H$. Then 
$$ \# (X'\backs H)=\# (X\backs H) + \# (B_0\backs H)\geq \#(X\backs H)+1
=\#(X\backs \span Y) +1, $$
where $Y\eqbd X\cap H$. Applying Corollary~\ref{cor_hom_i+}, we see that
$p_\eta^{\#(X'\backs H)}\in \calI_+(X)$. Thus all generators of $\calI(X')$
are in $\calI_+(X)$ and therefore $\calI(X')\subset \calI_+(X)$. Consequently,
$$\ker \calI_+(X) \subset \ker \calI(X') = \calP(X').$$

\nt
{\bf Step II.\/} We order $X'$ in a way that $B_0$ is placed after $X$,
we let $\{Q_B\eqbd p_{X'(B)} : B\in \B(X')  \}$ be the homogeneous basis 
for $\calP(X')$ (per the given order).  We define
$$\B' \eqbd  \B(X') \backs \B_+(X),\quad
F \eqbd  \span \Set{Q_B \,:\, B\in \B'}.$$
We will now prove that $F\cap \ker \calI_+(X) = \Set{0}$, hence that 
the quotient map $\calP(X')\to\calP(X')/F$ is an injection on $\ker \calI_+(X)$, and 
$$\dim\ker \calI_+(X)\le \dim\calP(X')-\dim F=\#\B_+(X)=\#\I(X).$$
Let 
$$f \eqbd  \sum_{B\in \B'} a(B) Q_B \in F. $$
Assume $f\in \ker \calI_+(X)$. 
We claim that $a(B)=0, \quad \forall B \in \B'$. To this end, we grade 
the bases in $\B'$ according to the location in $B_0$ of their maximal 
element, with respect to our fixed order. Note that the maximal element 
must be in $B_0$,  since otherwise $B\in\B(X)$. 
Assume that there exists $B_1 \in \B'$ such that $a(B_1) \neq 0$. Assume
further, without loss of generality, that $a(B)=0$ for every basis $B\in \B'$ with
higher grade. We then choose $H\eqbd \span(B_1\backs \max\{ B_1\} )$, let
$\eta$ be normal to $H$, and consider the differential operator
$D_\eta^k$, with 
$$k\eqbd \#\{x\in X'\backs H: x\not\in\span\{b\in B_1:\ b\preceq x\}\}.$$
Recall that the basis $B_1$ is {\it not\/} obtained by a greedy completion
of an independent subset of $X$ and that $B_0$ is in general position
with respect to $X$. This implies
$$ k \geq \#(X\backs Y)+1, \qquad {\rm where} \quad Y\eqbd X\cap H. $$
By Proposition~\ref{prop_hom_total},  $p_\eta^k  \in 
\calI_+(X)$, therefore  $D_\eta^k$ annihilates $\ker \calI_+(X)$ 
and, in particular, $f$. 

Finally, consider the set $\B''\subset\B'$ of bases  $B$ such that 
(i) their grade does not exceed the grade of $B_1$,  and
(ii) $D_\eta^k Q_B\not=0$. These are the bases with the property
$\# (X'(B)\backs H) \geq k$ or, equivalently, 
$$ \# \{ x \in X'\backs H : x\notin \spam 
\{b\in B : b\preceq x \} \}  \geq 
\# \{ x \in X'\backs H : x\notin \spam \{ b\in B_1 : 
b \preceq x \} \}. $$
Since no element of such a basis $B$ is located further than the
maximum  element $\max\{ B_1\}$ of $B_1$, each such basis $B$ must consist 
of a basis for $H$ augmented by the vector $\max\{ B_1\}$ itself. 
Hence, the set
$$\{D_\eta^k Q_B: B\in \B''\}$$
consists of (nonzero multiples of) elements of the homogeneous
basis of $\calP(X'\cap H)$. This implies that $a(B)=0$ for each $B\in\B''$, 
which leads to a contradiction, since $B_1\in\B''$.    \eop

We now state formally  the main theorem of this section.

\begin{theorem}\label{thm_main+}\hfill
\bd
\item{(1)}\quad 
$\dim \calP_+(X) = \dim \calD_+(X) = \#\I(X)$.
\item{(2)}\quad
The map $p \mapsto \inpro{p, \cdot}$ is a bijection between
$\calD_+(X)$ and $\calP_+(X)'$. 
\item{(3)}\quad
$\calD_+(X) =\Pi(V_+(X))= \ker \calJ_+(X)$. 
\item{(4)}\quad
The set $V_+(X,\lam)$ is correct for the space
$\calD_+(X)$, as well as for the space $\calP_+(X)$.\footnote{Note
that the set $V_+(X,\lam)$ depends on $X$, on $\lam$ and
on the augmented order basis $B_0$. The space $\calD_+(X)$
depends on $X$ and $B_0$, but not on $\lam$. Finally,
$\calP_+(X)$ depends only on $X$.}
\item{(5)}\quad
$\calP_+(X) = \ker \calI_+(X)$.
\item{(6)}\quad
$\calP_+(X)\oplus\calJ_+(X)=\Pi$.
\ed
\end{theorem}

\nt{\bf Proof.\/} The proof is analogous to that of Theorem~\ref{thm_main}.
 We put together inequalities obtained 
in Corollaries~\ref{cor_lb_d+}, \ref{cor_lb_p+}, \ref{cor_ub_p+}
and in Theorem~\ref{thm_dim_i+} to get
$$ \# \B_+(X)\leq \dim \calD_+(X) \leq \dim \calP_+(X) \leq \dim \ker\calI_+(X) 
\leq \# \B_+(X).   $$ 
This shows that equalities must hold throughout.
We then invoke Theorems~\ref{thm_pi_d+},~\ref{thm_sum+} and~\ref{thm_p_i+},
along with  Results~\ref{res_pi(sig)} and~\ref{res_sum} and 
Corollary \ref{cor_duala},  to prove  the remaining  claims of this 
theorem. \eop

\begin{theorem} \label{thm_intpoints+}
 Let $\calZ_+(X)$ be the integer points in the closed
zonotope $Z(X)$. Then $$\Pi(Z_+(X))=\calP_+(X)=\ker \calI_+(X),$$
provided that $X$ is unimodular.
\end{theorem} 

\nt{\bf Proof.\/} We first recall that, according to 
Result~\ref{thm_count}, 
$$\#\calZ_+(X)=\#\I(X)$$
in case $X$ is unimodular. That implies,
by invoking Theorem~\ref{thm_main+}, that
$$\dim\Pi(\calZ_+(X))=\#\B_+(X)=\dim\calP_+(X).$$
Hence our claim follows from the fact that
$$\Pi(\calZ_+(X))\subset \ker \calI_+(X).$$
The proof of this latter inclusion follows closely
that of Theorem~\ref{thm_intpoints},  hence is merely outlined: we need to
show that,
given any generator $q\eqbd p_{\eta_H}^{m(H)+1}$, $H\in\calF(X)$,
of $\calI_+(X)$, there exists $p\in \Pi$ that vanishes on $\calZ_+(X)$
and satisfies $p\most=q$. The existence of such $p$ follows from the
fact that, whatever facet hyperplane $H$ we choose, the set $\calZ_+(X)$
lies in the union
$$\cup_{j=0}^{m(H)}(a_j+H),$$
with $a_j\eqbd \sum_{k=1}^j y_k$, and where $\{y_j\}_{j=1}^{m(H)} \eqbd X\backs H$.
It is the unimodularity that guarantees that the above hyperplanes do not
depend on the order, and that the entire set $\calZ_+(X)$ lies in their
union.  (The description above assumes that $X\backs H$
all lie on the same side of $H$; the modifications that are needed
for the general case are notational.)  \eop

\subsection{\label{sec_basis+}Homogeneous basis and Hilbert series
for $\calP_+(X)$}

As before, we order $X$, and define, for each $I\in \I(X)$, 
$$X(I) \eqbd  \Set{x\in X \,:\, x \notin 
\span \Set{b\in I \,:\, b \le x}}.$$
Our goal now is to show that the set $$\Set{Q_I\eqbd p_{X(I)} \,:\, I \in \I(X)}$$
is a basis for $\calP_+(X)$.

\begin{theorem} \label{thm_basis+} \quad The set
$\Set{ Q_I\eqbd p_{X(I)} \,:\, I \in \I(X)}$ is a  basis for $\calP_+(X)$.
\end{theorem}

\nt{\bf Proof.\/}  Since the cardinality of the given set of polynomials
is $\#I(X)=\#B_+(X)=\dim\calP_+(X)$, and since obviously each one
of these polynomials lies in $\calP_+(X)$, 
we only need to show that the set $\{ Q_I : I\in \I(X)\}$ is linearly 
independent. For the proof of this part,
we order the augmented set $X'\eqbd  X\cup B_0$ such that $X$ retains
its internal order, and $B_0$ is placed after $X$.
Recall that each $I\in \I(X)$ has a well-defined extension
to a basis $\ex(I)\in\B(X')$:
$$\I(X) \to \B_+(X)\subset\B(X') \,:\, I \mapsto \ex(I).$$
We therefore examine the known homogeneous basis for $\calP(X')$.
The polynomials in that latter basis are $p_{X'(B)}$,
$B\in \B(X')$, with
$$X'(B)\eqbd \{x\in X': x\not\in\span\{b\in B: b\preceq x\}\}.$$
Now, given $I\in \I(X)$, since $\ex(I)\in \B(X')$ is a {\it greedy\/} 
extension of $I$, it easily follows that
$$X(I)=X'(\ex(I)),$$
hence that
$$Q_I\eqbd p_{X(I)}=p_{X'(\ex(I))}.$$
Thus, the set $\{Q_I\}_{I\in \I(X)}$ is a subset of the basis
$\{p_{X'(B)}: B\in\B(X')\}$ for $\calP(X')$, and the requisite linear
independence thus follows.
\eop

Note that the basis we just constructed is a homogeneous  extension of 
the homogeneous basis that was constructed for $\calP(X)$ in 
Section~\ref{sec_basis}. Moreover, the valuation 
$\val$ that was defined there on $\B(X)$ has just been extended 
also in the most natural way to $\I(X)$:
$$\val(I)\eqbd \#X(I),\quad I\in \I(X).$$
This motivates us to associate $X$ with 
{\it an external Hilbert series:\/}
$$h_+\eqbd h_{X,+}:\N\to \N\quad :\quad   j\mapsto \#\{I\in\I(X):\ \val(I)=j\}.$$
$h_+(j)$ equals thus to the dimension of $\calP_+(X)\cap \Pi_j^0$
and, by duality with $\calD_+(X)$, also to the dimension
of $\calD_+(X)\cap \Pi_j^0$.
The external Hilbert function is very special, in the sense
that its last non-zero entry always equals one. This fact is not easy
to observe by examining either of $\calI_+(X)$, $\calJ_+(X)$ or $\calD_+(X)$.
However, it trivially follows from the structure of $\calP_+(X)$: the
unique polynomial of maximal degree of the form $p_Y$, $Y\subset X$,
is $p_X$. One can use the above construction of a basis for $\calP_+(X)$
to conclude that $h_{X,+}(\#X-1)$ is the number of equivalence classes of
$X$ under the equivalence ($x\sim y$ iff $\{x,y\}$ is a dependent set).

\begin{example} Let
$$X=[x_1,x_2,x_3]\eqbd \left[\begin{array}{rrr}1&0&1\\
                                    0&1&-1\\
				    \end{array}\right].$$
This $X$ corresponds to a complete graph of three vertices,
and is unimodular, as is every graphical $X$. The zonotope $Z(X)$
has $7$ vertices in its closure. A basis for $\calP_+(X)$ is given 
by $$\{p_{Y}: Y\in 2^X\bks \{\{x_3\}\}\}.$$
The external Hilbert series is $h_{X,+}=(1,2,3,1)$.
With $x_4\eqbd x_3{\perp}=(1,1)'$, the ideal $\calI_+(X)$ is generated 
by the three polynomials $$p_{x_i}^3,\quad i=1,2,4.$$
It is clear that $h_{X,+}$ is indeed the correct Hilbert series
for this ideal. 

In contrast, the space $\calD_+(X)$ and its ideal $\calJ_+(X)$ are not
unique,
and depend on the choice of the augmented basis $B_0$. If we choose
$B_0=(y,z)$ with $y,z$ in general position with respect to $X$, then
the generators of $\calJ_+(X)$ become
$$p_{X\cup z}, \quad \hbox{and }\quad p_{(X\bks x)\cup y}, \ x\in X.$$
\eop
\end{example}

\medskip
Theorem \ref{thm_main+}\ yields the following characterization
$\calP_+(X)$: 

\begin{theorem}\label{thm_char+}\hfill
$$\calP_+(X)=\bigcap_{B\in\B(X)}\calP(X\cup B).$$
\end{theorem}

\nt{\bf Proof.\/}  The fact that $\calP_+(X)\subset \calP(X\cup B)$,
for any fixed basis $B$ for $\Rn$ follows trivially from the definitions of
$\calP(X)$ and $\calP_+(X)$ (and was used multiple times). We prove
that every polynomial $p$ in the intersection lies
in $\ker\calI_+(X)$ ($=\calP_+(X)$). Let $p$ be such a polynomial,
and let $H\in \calF(X)$. We need to show that, with $\eta_H \perp H$,
$D_{\eta_H}^{m(H)+1}p=0$, with $m(H)\eqbd \#(X\bks H)$. To this end,
we choose a basis $B\in\B(X)$ such that $\span(B\cap H)=H$ (such a basis
exists, since $H$ is a facet hyperplane.) Then,
with $X'\eqbd X\cup B$, $p\in \calP(X')=\ker \calI(X')$, hence is annihilated
by $D_{\eta_H}^{\#(X'\bks H)}$. Since only one vector of $B$ lies
outside $H$, we get that $\#(X'\bks H)=m(H)+1$, and the result follows.
\eop

\section{\label{sec_int}Internal algebra}

\subsection{\label{sec_main-}Main results}

We first recall the definition of internal bases.
We impose an (arbitrary but fixed) ordering $\prec$ on $X$. 
Let $B \in \B(X)$. If, for each $b \in B$, 
$$b \neq \max \Set{X\backs H}, \quad H \eqbd  \span \Set{B\backs
b}\in\calF(X),$$ 
then $B$ is called an {\it internal basis.\/} We denote the set of 
all internal bases by  $\B_-(X).$

For each $B \in \B(X)$, we define the {\it dual valuation\/} as follows: 
$$\val^*(B) \eqbd  \# \Set{b \in B \,:\, b \not = \max \Set{X \backs \span 
(B\backs b)}}.$$
Then,
$$\B_-(X) = \Set{ B\in \B(X) \,:\, \val^*(B)=n}.$$
We remind the reader about the fact mentioned in Section~\ref{sec_hyparr}:
$$ \# \B_-(X) = \sum_{I \in \I(X)} (-1)^{n-\#I}.$$

\nt{\bf Remarks.} (i) In matroid theory, $b$ is said to be {\it
internally active in\/} $B$ in the situation encountered above,  
i.e., if, for $b\in B$ and  $H\eqbd \span\{B\backs b\}$, we have that  
$b = \max \{ X\backs H \}$. The number $n-\val^\ast(B)$ is
known as the {\it internal activity\/} of $B$. (ii) The valuation $\val^\ast$
is then (matroid-)dual to the valuation $\val$ in the sense that it coincides
 with  the valuation $\val$ on the dual matroid of $X$. We make no use of this
duality since we do not develop zonotopal spaces on the dual matroid within 
this paper.

With a given order on $X$, we define the set of {\it barely
long subsets of\/} $X$:
$$L_-(X) \eqbd  \Set{ Y\subset X \,:\, Y\cap B \neq \emptyset, \quad B \in
\B_-(X)}.$$ 
The corresponding ideal is defined as
$$\calJ_-(X) \eqbd  \Ideal \Set{p_Y \,:\, Y \in L_-(X)},$$
and the notation for its kernel is
$$\calD_-(X) \eqbd  \ker \calJ_-(X).$$
It is clear that $\calJ_-(X) \supset \calJ(X)$ hence that $\calD_-(X) \subset
\calD(X)$. 
As is the case with external theory, one easily finds that
$\calD_-(X)$ depends on the ordering of $X$. (To recall, 
$\calD(X)$ does not depend on any ordering.)

Given a generic $X$-hyperplane arrangement $\calH(X,\lam)$,
we pick those vertices of it that are
associated with the internal bases $B\in \B_-(X)$
and call the resulting set $V_-(X,\lam)$. Then Theorem~\ref{thm_db'} and
its Corollary~\ref{cor_db'} apply to the space $\calD_-(X)$ and its
associated vertex set $V_-(X,\lam)$ to yield the following two results.

\begin{theorem} \label{thm_pi_d-}  $\Pi(V_-(X,\lam)) \subseteq \calD_-(X)$.  \eop
\end{theorem}

\begin{corollary} \label{cor_lb_d-}  $\dim \calD_-(X) \geq \# \B_-(X)$.  \eop
\end{corollary}

We now define a polynomial space and its ideal that will serve as duals
to the space  $\calD_-(X)$ and its ideal. 
\begin{eqnarray} 
\calP_-(X) & \eqbd & \cap_{x\in X}\calP(X\bks x),  \nonumber \\ 
\calI_-(X) & \eqbd & \Ideal \Set{p_{\eta_H}^{m(H)-1} \,:\, \eta_H \perp H, \quad 
H\in \calF(X)}. \label{defiminus}\end{eqnarray}
Note that $\calP_-(X)$ and $ \# \B_-(X)$ are independent of the order
$\prec$. Needless to say, the set  $\B_-(X)$ itself depends 
on that order.

\begin{theorem}\label{thm_char-}\hfill
$$\calP_-(X)=\ker\calI_-(X).$$
Moreover, for every $B\in \B(X)$,
$$\calP_-(X)=\cap_{x\in B}\calP(X\bks x).$$
\end{theorem}

\nt{\bf Proof.\/}  We prove this result by examining the corresponding
ideals: since $\calP(X\bks x)=\ker\calI(X\bks x)$, Theorem \ref{thm_main}, 
the stated result is proved once we show that (i) $\calI(X\bks x)\subset 
\calI_-(X)$, for every $x\in X$, and (ii) Given any $B\in \B(X)$,
$\calI_-(X)\subset\Ideal\{\cup_{b\in B}\calI(X\bks b)\}$.

For the proof of (i), fix $x\in X$, and denote $X'\eqbd X\bks x$. A generator $Q$
in the ideal $\calI(X')$ is of the form $Q\eqbd p_{\eta_H}^{\#(X'\bks H)}$, with
$H\in\calF(X')$. Then $H$ is also a facet hyperplane of $X$, and
obviously $\#(X'\bks H)\ge \#(X\bks H)-1$, and hence $Q\in \calI_-(X)$,
and (i) follows.

For (ii), we fix $B\in\B(X)$, and pick a generator of $\calI_-(X)$:
$Q=p_{\eta_H}^{m(H)-1}$, $H\in\calF(X)$, (\ref{defiminus}). We then 
choose $x\in B\bks H$, and denote $X'\eqbd X\bks x$.
Since $x\not\in H$, it is clear that $H\in\calF(X')$, and
then $\#(X'\bks H)=\#(X\bks H)-1=m(H)-1$. Thus the polynomial $Q$
lies in $\calI(X')$, and (ii) follows.
\eop

\medskip
Going back to the order we impose on $X$ (which is required for
the definition of $\calJ_-(X)$), we recall
the homogeneous construction of
a basis $(Q_B)_{B\in \B(X)}$ for $\calP(X)$ (see 
Theorem~\ref{thm_basis}), per that order.  Since $\B(X)$ is 
decomposed into internal and non-internal bases, it makes 
sense to decompose $\calP(X)$ accordingly, viz.,
$$\calQ_{in}(X)\eqbd \span\{Q_B: B\in \B_-(X)\},\quad \calQ_{ex}(X)\eqbd \span\{
Q_B:\ B\in \B(X)\bks\B_-(X)\}.$$
Then
$$\calP(X)=\calQ_{in}(X)\oplus\calQ_{ex}(X),$$
with the ``internal summand'' $\calQ_{in}(X)$ having the ``right
dimension'', i.e., $\#\B_-(X)$. However, in general that space differs
from $\calP_-(X)$ (whose dimension will be shown to equal $\#\B_-(X)$, too).
On the other hand, the complementary inclusion is true:

\begin{lemma}\label{lem_qex}
$\calQ_{ex}(X)\subset \calJ_-(X)$, hence
$$\codim J_-(X)\le \#\B_-(X).$$
\end{lemma}

\nt {\bf Proof.\/} Fix $B\in\B(X)\bks \B_-(X)$. Then $B$ contains
an internally active $b$: with $H\eqbd \span(B\bks b)\in\calF(X)$,
$b$ is the last vector in $X\bks H$. Examining the definition (\ref{defxb})
of the set $X(B)$,  it is then clear that $X\bks (X(B)\cup 
B)\subset H$.  (Indeed, if $x\in X\bks (H\cup B)$, then $x\prec b$, hence 
it is impossible that $x\in\span\{b'\in B: b'\prec x\}$, since the
latter span lies in $H$.) Next, it easily follows that $b$ belongs to,
and is internally active in every basis
$B'\subset X\bks X(B)$. Thus $X\bks X(B)$ does not contain an internal
basis, hence $X(B)\in L_-(X)$.

Recall now that $\Pi=\calP(X)\oplus\calJ(X)$ according to 
Theorem~\ref{thm_main}, and 
since $\calJ(X)\subset \calJ_-(X)$, we conclude that
$\calQ_{ex}(X)+\calJ(X)\subset \calJ_-(X)$, hence
$$\Pi=\calQ_{in}(X)\oplus\calQ_{ex}(X)\oplus\calJ(X)=\calQ_{in}+\calJ_-(X).$$
Consequently, 
$$\codim \calJ_-(X)\le \dim\calQ_{in}(X)=\#\B_-(X).$$
\eop

\begin{theorem}\label{thm_dimd-}
$$\dim \calD_-(X)= \#\B_-(X),$$
and  $\Pi(V_-(X,\lam))=\calD_-(X)$.
\end{theorem}

\nt {\bf Proof.\/}  Since $\calD_-(X)=\ker \calJ_-(X)$, by definition,
we have that $\dim\calD_-(X)=\codim\calJ_-(X)$.
However, by Corollary \ref{cor_lb_d-}, $\dim\calD_-(X)\ge \#\B_-(X)$, while 
by Lemma~\ref{lem_qex}, 
$\codim\calJ_-(X)\le \#\B_-(X)$. Thus:
$$\dim\calD_-(X)=\#\B_-(X).$$
But then, $\Pi(V_-(X,\lam))$  is a subspace of $\calD_-(X)$  
(Theorem \ref{thm_pi_d-}) 
of dimension $\#\B_-(X)$, so we have that
$\Pi(V_-(X,\lam))=\calD_-(X)$.
\eop

\begin{corollary}
$\calJ_-(X)=\calJ(X)\oplus \calQ_{ex}(X).$
\end{corollary}

\begin{theorem} \label{thm_sum-}  $\calP_-(X)+\calJ_-(X)=\Pi$. 
In particular, $\dim\calP_-(X)\ge \#\B_-(X)$.
\end{theorem}

\nt {\bf Proof.\/}
Since we already know that $\codim\calJ_-(X)=\#\B_-(X)$, the second
claim in the theorem follows from the first. Let us thus prove the
first.

From Lemma~\ref{lem_qex}, we know that $\calQ_{in}(X)+\calJ_-(X)=\Pi$. 
Thus, it is enough to show that $\calQ_{in}\subset \calP_-(X)+\calJ_-(X)$.
We achieve this latter relation by showing that every polynomial
$Q_B$, $B\in\B_-(X)$, lies in $\calP_-(X)+\calJ_-(X)$, and use the
following general approach. Fixing $B\in \B_-(X)$, we know that
$Q_B=p_{X(B)}$, for  suitable $X(B)\subset X$. We decompose
$X(B)$ in a certain way $X(B)=Z\cup W$. Thus
$$Q_B=p_Zp_W.$$
We then replace each $w\in W$ by a vector $w'$ (not necessarily
in $X$), to obtain a new polynomial
$$\til{Q_B}\eqbd p_{Z}p_{{W'}},$$
and prove that (i) $\til{Q_B}\in \calP_-(X)$, and (ii)
$Q_B-\til{Q_B}\in \calJ_-(X)$.

So, let $Q_B=p_{X(B)}$ be given. If $Q_B\in\ker \calI_-(X)=\calP_-(X)$,
there is nothing to prove. Otherwise, let ${\bf H}\subset\calF(X)$ be 
the collection of {\it all\/} facet hyperplanes for which
$D_{\eta_H}^{m(H)-1}Q_B\not=0$. The set ${\bf H}$ is not empty, since otherwise
$Q_B\in \ker \calI_-(X)$. Given $H\in {\bf H}$, we conclude
that $\#(X(B)\bks H)\ge m(H)-1$, hence that, with $Y\eqbd X\bks X(B)$,
$\#(Y\bks H)\le 1$. Since $B\subset Y$, the set
$Y\bks H$ must be a singleton $x_H\in B$. We denote
$$X_\bfH\eqbd \{x_H:\ H\in \bfH\}.$$

Define
$$W\eqbd \{\max \{ X\bks H \} :\ H\in\bfH\}.$$
We index the vectors in $W$ according to their order in $X$:
$W=\{w_1\prec w_2\prec\ldots\prec w_k\}.$
For each $1\le i\le k$, we define 
$$X_i\eqbd \{x_H:\ H\in\bfH,\ \max\{X\bks H\}= w_i\},\quad \bfH_i\eqbd \{H\in
\bfH:\ x_H\in X_i\}.$$
Thus, $X_\bfH=\bigcup_{i=1}^kX_i$.

Setting all these notations, we first observe that $W\cap X_\bfH=\emptyset$,
i.e., $w_i$ does not lie in $X_i$. Indeed, the set $X_\bfH$ is a subset
of every $B'\in \B(Y)$, with $\span(B'\bks x_H)=H$ for each $x_H\in X_\bfH$.
If some $x_H$ is $\max\{ X\bks H\}$, it will be internally active in every
$B'\in\B(Y)$, which would imply that $\B(Y)$ does not contain internal bases,
which is impossible since $B\in \B(Y)$. Thus, $W\subset X(B)$, and we
define $Z\eqbd X(B)\bks W$, to obtain
$$Q_B=p_Zp_W.$$

Define further:
$$S_i\eqbd \cap\{H: H\in\cup_{j=1}^i \bfH_j\},\quad S_0\eqbd \Rn.$$
Then, for $i=1,\ldots,k$,
$S_{i-1}=S_{i}\oplus \span\, X_{i}$, 
and $w_{i}\in S_{i-1}\bks S_{i}$.  
Thus, for $i=1,\ldots k$, the vector $w_i$ admits a unique representation of the
form
\begin{equation}\label{defxijprime}
w_i=w_i'+\sum_{x\in X_{i}}a_xx,\quad w_i'\in S_{i},\  a_x\in \R\bks
\{0\}.
\end{equation}
Define
$$W'=\{w_1',\ldots,w_k'\},\ \hbox{and }
\tilQ_B\eqbd p_{Z}p_{W'}.$$
We prove first that
$$\tilQ_B-Q_B=p_{Z}(p_{W'}-p_W)$$
lies in $\calJ_-(X)$. 
To this end, we multiply out the product
\begin{equation}
p_{W'}=\prod_{i=1}^k p_{w_i'}=\prod_{i=1}^k(p_{w_i}-\sum_{x\in X_i
}a_x p_{x}).\end{equation} 
Every summand in the above expansion is of the form $p_\Xi$, with
$\Xi$ a suitable mix of $W$-vectors and $X_\bfH$-vectors.
The summand $p_W$ in the above expansion in canceled when we subtract
$Q_B$. Any other $\Xi$ is obtained from $W$ by replacing at least once 
a  $w_i$ vector by some vector  in $X_i$, which we denote by $x_i$. Let
$w_{i_1}\prec w_{i_2}\prec\ldots\prec w_{i_j}$ be all the $w$-vectors
in $W\bks \Xi$, and let $H_1$ be the facet hyperplane that
corresponds to $x_{i_1}$ ($H_1\eqbd \span(B\bks x_{i_1})$.)
Then, we have that $w_{i_1}\in X\bks(Z\cup \Xi)=:Y'$, and we claim that
$Y'\bks w_{i_1}\subset H_1$. To this end, we write
$Y'\bks H_1=((Y'\cap Y)\bks H_1)\cup (Y'\bks Y)\bks H_1$.
Now, $Y\bks H_1=x_{i_1}$, and since $x_{i_1}\not\in Y'$ (as it was replaced
by $w_{i_1}$), the term $(Y'\cap Y)\bks H_1$ is empty. The second
term consists of $(w_{i_m})_{m=1}^j\bks H_1$. However,
$w_{i_m}\in S_{i_m-1}\subset S_{i_1}\subset H_1$, for every $m\ge 2$.
Thus, $w_{i_1}$ is the only vector in $Y'\bks H_1$. Being also the
last vector in $X\bks H_1$, we conclude that $w_{i_1}$ is
internally active in every $B\in \B(Y')$, hence
that $p_{Z\cup W}\in \calJ_-(X)$. This being true for every summand
in $\tilQ_B-Q_B$, we conclude that this latter polynomial lies in
$\calJ_-(X)$.

We  now prove that $\tilQ_B=p_{Z\cup W'}\in \ker \calI_-(X)$. To this end,
we need to show that, for every $H\in\calF(X)$, 
$\#((Z\cup W')\bks H)<m(H)-1$. We divide the discussion here to three cases.
As before, $Y\eqbd X\bks X(B)$.

Assume first that $H\in \bfH_i$ for some
$1\le i\le k$. Then, for $X(B)=Z\cup  W$ we had that
$\#((Z\cup W)\bks H)=m(H)-1$. Now, $x_H$ is the only
vector in $Y\bks H$, and $x_H\in X_i$.  Thus,
the subset $X_j\subset Y$, must lie in $H$ for every $j\not= i$,
which means that we conclude that,
$w_j\in H$ iff $w'_j\in H$ (since
$w_j-w_j'\in \span X_j\subset H$).
Finally, while $w_i\not\in H$, $w'_i\in S_i\subset H$,
hence, altogether, $\#(W'\bks H)<\#(W\bks H)$,  and we
reach the final conclusion that
$$\#((Z\cup W')\bks H)<
\#((Z\cup W)\bks H)=m(H)-1.$$

Secondly, assume that $S'\eqbd S_k\cap H\not= S_k$. Then, necessarily,
$U\eqbd (Y\cap S)\bks S'$ contains at least two vectors (otherwise, 
all the vectors of $Y$ but one lie in the rank deficient
set $(Y\cap S')\cup (Y\bks S_k)$). Now, with
$$m_1\eqbd \#\{w\in W: w\in H\wedge w'\not \in H\}, \hbox{ and } 
m_2\eqbd \#(\cup_{i=1}^k(X_i\bks H)),$$ 
we know that $\#((Z\cup W')\bks H)=m(H)-\#U+m_1-m_2$. However, we must
have that $m_1\le m_2$: if
$w_i'\not\in H$, while $w_i\in H$, then, since 
$w_i-w_i'\in \span X_i$, we have that $\#(X_i\bks H)>0$.

Finally, we assume $H\in\calF(X)\bks \bfH$, and $S_k\subset H$.
Let $j\ge 1$ be the minimal index $i$ for which $S_i\subset H$. We modify
the definition of $m_1$ and $m_2$ from the second case by replacing
$W$ by $W\bks w_j$ in the definition of $m_1$, and removing $X_j$
from $\cup_{i=1}^k X_i$ in the definition of $m_2$. We still
have that $m_1\le m_2$, by the same argument as above. However, the
set $U$ that we used in the previous case is not available for us.
Instead, we examine the relation
$$w_j-w'_j\in\span X_j.$$
We know {\it a priori\/} that $S_j\oplus \span X_j=S_{j-1}$.
Since $S_j\subset H$, while $S_{j-1}\not\subset H$,
we must have that $X_j\bks H\not=\emptyset$. But, $w_j'\in H$,
hence, with $U\eqbd (w_j\cup X_j)\bks H$, $\#U\ge 2$, hence the argument
used for the previous case works here, with $U$, $m_1$ and $m_2$
modified as explained.
\eop

We now establish our second non-trivial theorem in this section.

\begin{theorem} \label{thm_dim_i-}
$$\calP_-(X)\cap\calQ_{ex}(X)=\{0\}.$$
In particular, $\dim\calP_-(X)\le \#\B_-(X)$.
\end{theorem}

\nt {\bf Proof.\/}  The second  claim follows from the first:
since $\calP_-(X)\subset \calP(X)=\calQ_{in}(X)\oplus \calQ_{ex}(X)$,
the first claim implies that 
$$\dim\calP_-(X)\le \dim \calP(X)/\calQ_{ex}(X)=\dim\calQ_{in}(X)=\#\B_-(X).$$

In order to prove the first claim, we denote
$$\oB'\eqbd \oB(X)\bks\oB_-(X).$$
Then, with $Q_B$ the polynomial in the homogeneous basis
for $\calP(X)$ that corresponds to  $B\in\oB(X)$, 
we pick a generic function $f\in \calQ_{ex}(X)$:
$$f=\sum_{B\in\oB'}a(B)Q_B,$$
and assume that $f\in\ker \calI_-(X).$ 
The proof that $f=0$ will be done as follows.  In addition
to the existing order $\prec$ on $X$, we will impose a
full order $\bprec$ on the bases in $\oB'$. Assuming  $f\neq 0$,
we will then select $B'\in \oB'$ which is minimal (in the full
order $\bprec$) with respect to the condition $a(B)\neq 0$. Thus
\begin{equation}
f-a(B')Q_{B'}=\sum_{B\in\oB',\ B\subbsucc B'}a(B)Q_B.\label{temp}
\end{equation}
We will then select a facet hyperplane $H\in \calF(X)$,  
and, with $\eta\in \Rn$
a normal to that hyperplane, apply to both
side of (\ref{temp}) the differential operator $D_\eta^{m(H)-1}$. Since
$f\in \ker\calI_-(X)$, by assumption,
$D_\eta^{m(H)-1}f=0$. Therefore,
$$-D_\eta^{m(H)-1}a(B')Q_{B'}=
\sum_{B\in\oB',\ B\subbsucc B'}a(B) D_\eta^{m(H)-1} Q_B.$$
The key of the proof will be to show, with the proper selection
of the full order, and with proper selection of the hyperplane
$H$, that the polynomial $D_\eta^{m(H)-1}Q_{B'}$ is independent of the
polynomials $D_\eta^{m(H)-1} Q_B$, $B\in\oB'$, $B\bsucc B'$.
This will imply that $a(B')=0$, hence will provide the sought-for
contradiction.

Here are the details.
We start with the introduction of the order on $\oB'$. To this end,
given $B\in \oB'$, and $b\in B$, we recall that $b$ to is said to be
 {\it internally active in\/} $B$ if $b$ is the maximal vector in 
$X\bks \spa\{B\bks b\}$.
We denote by $\alp(B)$ the number of internally active vectors in $B$.
Note that $B\in \oB'$ if and only if $\alp(B)>0$. We choose the order
on $\oB'$ to respect the number of internally active vectors, i.e.,
$$B\bprec\tilB\implies \alp(B)\le \alp(\tilB).$$

Next, let $B\in \oB'$ and $x\in B$. Set $H\eqbd \spa(B\bks x)$.
We say that $x$ is an {\it $H$-shield\/}
of $B$ if $X\bks\{X(B)\cup x\}$ is not full rank (hence lies
in $H$), but $x$ is not the maximal vector in $X\bks H$. 

Now, in order to show that $a(B')=0$, we select any internally active
$b'\in B'$ (there must be at least one, since $B\in\oB'$), define
$H'\eqbd \spa(B'\bks b')$, and take $\eta$ to be normal to $H'$.
With $G\eqbd D^{m(H')-1}_{\eta}$, we know that $Gf=0$ (since
$f\in\ker\calI_-(X)$, by assumption). Moreover,
since $m(H')=\#(X\bks H')$, and since, for any $B\in\oB(X)$,
$X(B)$ is disjoint of $B$ (while $B$ contains at least one vector
from $X\bks H'$), we have $\#(X(B)\bks H')<m(H)$. If 
$\#(X(B)\bks H')<m(H)-1$, then $GQ_B=0$. Otherwise, up to a 
non-zero multiplicative constant, $GQ_B=p_{X(B)\cap H'}$.
Note that, with $X'\eqbd X\cap H'$,  we get $X(B)\cap H'=X'(B\cap H')$,
i.e., the set  $B\cap H'$  spans $H'$ (otherwise, $GQ_B=0$), hence 
lies in $\oB(X\cap H')$; with $X'$ retaining its $X$-order $\prec$, 
the construction of  a homogeneous basis for $\calP(X\cap H')$ associates 
the basis $B\cap H'$ with the polynomial $p_{X(B)\cap H'}$, i.e., with 
the polynomial $GQ_B$ (up to the aforementioned constant).  
The selection of $H'$ clearly implies that $GQ_{B'}\not=0$, hence, 
in particular, that $GQ_{B'}=c'p_{X(B')\cap H'}$ for some $c'\neq 0$.  
Set $$\oB''\eqbd \{B\in \oB': B\bsucc B',\ \#(X(B)\bks H')=m(H')-1\}.$$
We get
\begin{equation} 
-c'a(B')p_{X(B')\cap H'}=\sum_{B\in\oB''}c(B) a(B) p_{X(B)\cap
H'}.\label{tempa}
\end{equation}
By our argument above, all the polynomial summands on both
sides of~(\ref{tempa})
belong to the homogeneous basis of $\calP(X')$. However, {\it a priori\/}
we cannot conclude that all the coefficients in~(\ref{tempa})
equal $0$, since we have not excluded the possibility that polynomials
from the aforementioned homogeneous basis make multiple appearances
in~(\ref{tempa}). Since we only need to prove that $a(B')=0$,
we need only to show that $p_{X(B')\cap H'}$ is not one of the
summands in the right hand side of~(\ref{tempa}). This is equivalent
to proving that, for each $B\in \oB''$, $B'\cap H'\neq B\cap H'$. 

Let us assume that $B\in \oB''$ and $B'\cap H'=B\cap H' \bdeq A$.
Obviously, $b \eqbd B\bks A\neq  b'$ (otherwise, $B=B')$, hence
$b$ is an $H'$-shield of $B$. We will show that the existence
of $H'$-shield in $B$ implies that $\alp(B)<\alp(B')$,
which will contradict the assumption that $B'\bprec B$.

It remains to prove the crucial thing: that $\alp(B)<\alp(B')$.
The argument for that is as follows: We recall that
$B\bks b=B'\bks b'=:A$, and $A$ is a basis for $H'$. We already know
that $b$ is not internally active in $B$,
while $b'$ is internally active in $B'$: $b'$ is the last vector in
$X\bks H'$, and $b \prec b'$.
We will show that if $x\in A$ is internally active in $B$ then it 
is also internally active in $B'$. Then, all the internally active
vectors in $B$ are internally active in  $B'$, while $B'$ contains 
an additional internally active element, viz., $b'$.     

So, let $x\in A$ be internally active in $B$. Set $S\eqbd A\bks x$.
Note that $\rank S=n-2$. If $x$ is not internally active
in $B'$, then there exists $y \succ x$ such that $y\not\in 
\spa\{S\cup b'\}$. Assume $y$ to be maximal element outside $\spa\{S\cup 
b'\}$. We get the contradiction to the existence of such $y$ by
showing that it is impossible to have $y \succ b'$, and it is also
impossible to have $y \prec b'$.

If $y\succ b'$, then, since $b'$ is maximal outside $\spam\{B\bks b'\}
=\spam A=\spam\{S\cup x\}$, we have that $y\in\spam\{S\cup x\}$. Also,
since $y\succ x$, and $x$ is maximal outside $\spam \{ B\bks x\}=
\spam\{S\cup b\}$, we have $y\in\spam\{S\cup b\}$. But $S\cup b\cup x=B$, 
and $B$ is independent, hence $y\in \spam S$, which is impossible since 
we assume $y$ to be outside $\spam\{S\cup b'\}$.

Otherwise, $y\prec b'$. The maximality of $y$ then implies
that $x\prec y \prec b'$. The maximality of $x$ outside $\spam\{S\cup b\}$
implies that $b'\in\spam\{S\cup b\}$. Since $b'\not\in S$,
we obtain that $\spam\{S\cup b\}=\spam\{S\cup b'\}$, which is impossible
since $y$ lies in exactly one of these two spaces. 
\eop

We now state formally  the main theorem of this section.

\begin{theorem}\label{thm_main-}\hfill
\bd
\item{(1)}\quad 
$\dim \calP_-(X) = \dim \calD_-(X) =\#\B_-(X)$.
\item{(2)}\quad
The map $p \mapsto \inpro{p, \cdot}$ is a bijection between
$\calD_-(X)$ and $\calP_-(X)'$. 
\item{(3)}\quad
$\calD_-(X) =\Pi(V_-(X,\lam))= \ker \calJ_-(X)$. 
\item{(4)}\quad
The vertex set $V_-(X,\lam)$ is correct for $\calD_-(X)$ as well as for
$\calP_-(X)$.
\item{(5)}\quad
$\calP_-(X) = \ker \calI_-(X)$.
\item{(6)}\quad
$\calP_-(X)\oplus\calJ_-(X)=\Pi$.
\ed
\end{theorem}

\nt{\bf Proof.\/} This proof is not directly parallel but still
quite similar to that  of  Theorem~\ref{thm_main}.  We put together 
inequalities and equalities obtained in Corollary~\ref{cor_lb_d-} and 
Theorems~\ref{thm_char-}, \ref{thm_dimd-} and~\ref{thm_dim_i-} to get
$$ \# \B_-(X) = \dim \calD_-(X) \leq \dim \calP_-(X) = \dim 
\ker\calI_-(X) \leq \# \B_-(X).   $$ 
This shows that equalities must hold throughout. We then invoke 
Theorems~\ref{thm_char-},~\ref{thm_dimd-} and~\ref{thm_sum-}, along with 
Result~\ref{res_pi(sig)}, Corollary~\ref{cor_duala} and Result~\ref{res_sum},
  to obtain the remaining  claims of this theorem. \eop

\begin{theorem} \label{thm_intpoints-}
 Let $\calZ_-(X)$ be the integer points in the interior of the
zonotope $Z(X)$. Then $$\Pi(Z_-(X))=\calP_-(X)=\ker \calI_-(X),$$
provided that $X$ is unimodular.
\end{theorem} 

\nt{\bf Proof.\/} The proof is analogous to the proofs 
of Theorems~\ref{thm_intpoints} and~\ref{thm_intpoints+} before. 
We first recall the count  $$\#\calZ_-(X)=\#\B_-(X),$$
which is true for a unimodular $X$. That implies,
by invoking Theorem~\ref{thm_main-}, that
$$\dim\Pi(\calZ_-(X))=\dim\calP_-(X).$$
Hence our claim follows from the fact that
$$\Pi(\calZ_-(X))\subset \ker \calI_-(X).$$
The proof of this latter inclusion requires us to show  
that, given any generator $q\eqbd p_{\eta_H}^{m(H)-1}$, $H\in\calF(X)$,
of $\calI_-(X)$, there exists $p\in \Pi$ that vanishes on $\calZ_-(X)$
and satisfies $p\most=q$. The existence of such $p$ follows from the
fact that, whatever facet hyperplane $H$ we choose, the set $\calZ_-(X)$
lies in the union
$$\cup_{j=1}^{m(H)-1}(a_j+H),$$
with $a_j\eqbd \sum_{k=1}^j x_k$, and where $\{x_j\}_{j=1}^{m(H)}=X\backs H$;
the hyperplanes in the above union do not depend on the order we impose
on $X\backs H$. As before, we can assume without loss of generality that 
$X\backs H$ all lie on the same side of $H$.  \eop

\subsection{\label{sec_basis-}Homogeneous basis and Hilbert series for $\calP_-(X)$}

The {\it internal Hilbert series\/} $h_{X,-}$ records the homogeneous
dimensions of $\calP_-(X)$:
$$h_-(j)\eqbd h_{X,-}(j)\eqbd \dim(\calP_-(X)\cap\Pi_j^0)=
\dim(\calD_-(X)\cap\Pi_j^0), \quad j\in\N.$$
While it is not true in general that
the polynomials 
$Q_B\eqbd p_{X(B)}$, $B\in\B_-(X)$,
form a basis for $\calP_-(X)$, they {\it can\/} be used for computing
$h_{X,-}$:
$$h_{X,-}(j)=\#\{B\in \B_-(X): \val(B)=\#X(B)=\deg Q_B=j\}.$$
In other words, the homogeneous dimensions of the (order-dependent)
space $\calQ_{in}(X)$ coincide with those of $\calP_-(X)$:
$$\dim(\calQ_{in}(X)\cap \Pi_j^0)=\dim(\calP_-(X)\cap\Pi_j^0),\quad \forall j.$$
The simplest way to observe this fact, is to follow the proof of 
Theorem \ref{thm_sum-}: Every $Q_B$ there was proved to be writable as
$$Q_B=\tilQ_B+f_B$$
with $f_B\in \calJ_-(X)$ and $\tilQ_B\in\calP_-(X)$. The fact
that $\tilQ_B$, $B\in\oB_-(X)$, are independent follows directly
from the independence of $Q_B$, $B\in \oB_-(X)$, and the fact
that the sum $\calQ_{in}(X)+\calJ_-(X)$ is direct. Since we know by now that
$\dim\calP_-(X)=\#\B_-(X)$, we conclude that

\begin{corollary}
The polynomials $\tilQ_B$, $B\in \B_-(X)$, from the proof of
Theorem \ref{thm_sum-} form a basis for $\calP_-(X)$.
\end{corollary}

Now, each $\tilQ_B$ is obtained by replacing some of the factors
$p_w$, $w\in  X$ of $Q_B$, by polynomials $p_{w'}$, $w'\in\Rn\bks 0$.
Thus, trivially, $\deg Q_B=\deg\tilQ_B$, hence we may indeed compute
$h_{X,-}$ via the polynomials $(Q_B)_{B\in\B_-(X)}$.

The fact that the spaces $\calQ_{in}(X)$ and $\calP_-(X)$ are different
is somewhat less trivial. For example, in two dimensions they are actually
the same. In three dimensions, however, they may not be the same, as
the following example shows:

\begin{example} Let
$$X=[x_1,\ldots,x_5]\eqbd \left[\matrix{1&0&0&1&1\cr
                                    0&1&0&2&1\cr
				    0&0&1&1&1\cr}\right].$$
Then
$$\B_-(X)=\{[x_1,x_2,x_3],[x_1,x_3,x_4],[x_1,x_2,x_4]\}\bdeq (B_1,B_2,B_3).$$
Our theory asserts, then, that $\dim\calP_-(X)=3$. Indeed, one
verifies directly that 
$$\calP_-(X)=\span\{1,p_{x_2},p_{x_4}\}.$$
The polynomials $Q_{B_1}=1$, $Q_{B_2}=p_{x_2}$ and $Q_{B_3}=p_{x_3}$
span the space
$\calQ_{in}(X)=\span\{1,p_{x_2},p_{x_3}\}$. The two spaces,
$\calP_-(X)$ and $\calQ_{in}(X)$ are different, but they produce
the same Hilbert series: $$h_{X,-}=(1,2,0,0,\ldots).$$
\qed
\end{example}

\section{\label{sec_conclude}Concluding remarks}

A key component of the theory of zonotopal algebra is the explicit
use of the polynomials $p_Y$, $Y\subset X$, in the construction of
the $\calJ$-ideals, as well as in the construction of the
$\calP$-spaces. In the context of the $\calJ$-ideal, the only
deviation is the use of the external basis $B_0$ for defining
$\calJ_+(X)$: a generator of $\calJ_+(X)$ is of the form
$p_Yp_b$, with $Y\subset X$, and $b\in B_0$. For example,
if $\rank(X\bks Y)=n-1$, then $b$ is the first vector in 
$B_0\bks \span(X\bks Y)$.

In the context of the $\calP$-spaces, the deviation from a direct
use of polynomials $p_Y$ occurs in the case of $\calP_-(X)$.
In an earlier formulation of the internal theory, we defined
the internal $\calP$-space as
$$\til{\calP_-(X)}\eqbd \span\{p_Y: Y\in S_-(X)\},$$
with the subset of {\it very short $X$-sets\/} defined as
$$S_-(X)\eqbd \{Y\subset X: \rank(X\bks (Y\cup x))=n,\ \forall x\in X\bks Y\}.$$
While this variant {\it is\/}  spanned by polynomials of the form
$p_Y$, $Y\subset X$, and while it is straightforward to check
that this space is a subspace of $\calP_-(X)$, we did not prove
that the two variants coincide. We conjecture, however, that the 
two spaces do coincide:

\begin{conjecture}\label{pminusconj}
For every $X$, $\calP_-(X)=\span\{p_Y: Y\in S_-(X)\}.$
\end{conjecture}

Note that proving the above conjecture is tantamount to showing that
the polynomials $p_Y\in \calP_-(X)$, $Y\subset X$, span $\calP_-(X)$: 
a polynomial $p_Y$, $Y\subset X$ lies in $\calP_-(X)$ iff $Y$ is very 
short. In any event, the proof of Theorem~\ref{thm_sum-} reveals the 
following information about $\calP_-(X)$:

\begin{corollary}
The space $\calP_-(X)$ is spanned by polynomials of the form
$p_Y$, $Y\subset \Rn$. Moreover, there exists a basis for 
$\calP_-(X)$ such that each polynomial in that basis is of the form
$$p_Y p_Z,$$ with $Y\subset X$, $Z\subset \Rn$, and $\#Z\le n-2$.
\end{corollary}

\medskip
A second remark concerns the $\calI$-ideals and $\calD$-spaces.
While the $\calI$-ideals admit a simple set of generators,
we do not know of any simple algorithm for constructing an
{\it explicit\/} basis of a $\calD$-space. Another remark,
of a different flavor, concerns a special property of the 
central $\calD(X)$ space: $\calD(X)$ is the smallest 
translation-invariant subspace of $\Pi$ that contains
$\calD(X)\cap \Pi_{\#X-n}^0$. This property does not extend to 
other $\calD$-spaces. For example, with $N\eqbd \#X$, 
$\dim(\calD_+(X)\cap\Pi_N^0)=1$, while 
$\dim(\calD_+(X)\cap\Pi_{N-1}^0)>n$, unless $X$ is a tensor 
product, i.e., consists of $n$ different vectors, each appearing 
with arbitrary multiplicity.

\section*{Acknowledgments}

The authors are grateful to Carl de Boor, Nira Dyn, Uli Reif,
Frank Sottile and Bernd Sturmfels for fruitful discussions
and to the referees for helpful comments.

\bibliographystyle{acm}
\bibliography{zono}

\end{document}